\documentclass[9pt,lettersize,leqno,fullpage]{article}
\usepackage[usenames,dvipsnames]{color}
\usepackage{amsmath, amssymb, amscd, amsthm, amsfonts, amsbsy, appendix}
\usepackage{bbm, stmaryrd, mathrsfs, bm, linearb, upgreek, skak, textcomp, mathdots}
\usepackage{mathtools, titlesec, extarrows, tocloft, url, dsfont}
\usepackage[intoc]{nomencl}
\usepackage[nottoc]{tocbibind}
\usepackage[all]{xy} 
\usepackage{graphicx, pgf}
\usepackage[T1]{fontenc}

\DeclareMathAlphabet{\mathpzc}{OT1}{pzc}{m}{it} 

\theoremstyle{plain}
\newtheorem{thm}{Theorem}[section]
\newtheorem{cor}[thm]{Corollary}
\newtheorem{lmm}[thm]{Lemma}
\newtheorem{prpn}[thm]{Proposition}
\newtheorem{rem}[thm]{Remark}

\theoremstyle{definition}
\newtheorem{defn}[thm]{Definition}
\newtheorem{eg}[thm]{Example}

\newtheorem*{conj}{Conjecture}
\newtheorem*{thma}{Theorem A}
\newtheorem*{thmb}{Theorem B}
\newtheorem*{thmc}{Theorem C}

\newcommand{\CP}{\mathbb{CP}}

\numberwithin{equation}{section}

\newcommand{\Q}{\mathbb{Q}}
\newcommand{\R}{\mathbb{R}}

\newcommand{\Z}{\mathbb{Z}}

\newcommand{\into}{\hookrightarrow}

\newcommand{\bgd}{\begin{displaymath}}
\newcommand{\edd}{\end{displaymath}}
\newcommand{\bge}{\begin{equation}}
\newcommand{\ede}{\end{equation}}
\newcommand{\bgea}{\begin{eqnarray*}}
\newcommand{\edea}{\end{eqnarray*}}
\newcommand{\bgeA}{\begin{eqnarray}}
\newcommand{\edeA}{\end{eqnarray}}
\newcommand{\bgc}{\begin{center}}
\newcommand{\edc}{\end{center}}
\newcommand{\ben}{\begin{enumerate}}
\newcommand{\een}{\end{enumerate}}
\newcommand{\bgi}{\begin{itemize}}
\newcommand{\edi}{\end{itemize}}
\newcommand{\hf}{\hspace*{0.5cm}}

\newcommand{\lp}{\bullet}

\title{{\LARGE Homotopy groups and periodic geodesics of closed 4-manifolds}}
\author{Samik Basu, Somnath Basu}
\date{}
\usepackage{pdfpages}
\usepackage[pdfauthor={Samik Basu, Somnath Basu},%
        pdftitle={Homotopy groups and periodic geodesics of closed 4-manifolds},%
        pdftex]{hyperref}
\hypersetup{
	colorlinks,
	citecolor=red,
	filecolor=black,
	linkcolor=blue,
	urlcolor=black
}

\begin{document}

\maketitle

\begin{abstract}
Given a simply connected, closed four manifold, we associate to it a simply connected, closed, spin five manifold. This leads to several consequences : the stable and unstable homotopy groups of such a four manifold are determined by its second Betti number, and the ranks of the homotopy groups can be explicitly calculated. We show that for a generic metric on such a smooth four manifold with second Betti number at least three the number of geometrically distinct periodic geodesics of length at most $l$ grows exponentially as a function of $l$. The number of closed Reeb orbits of length at most $l$ on the spherization of $T^\ast M$ also grow exponentially for any Reeb flow. 
\end{abstract}

\vspace*{0.5cm}
\tableofcontents
\vspace*{0.75cm}


\section{Introduction}\label{intro}

\hf\hf Questions regarding classification of manifolds present themselves naturally in diverse examples. There are questions that, in disguise, straddle the worlds of topology and geometry. For each conjugacy class in the fundamental group of a Riemannian manifold, we get a periodic geodesic\footnote{In the literature the term {\it closed geodesics} is often used instead of periodic geodesics.} by looking at a closed curve representatives which minimize length. However, these may not be geometrically distinct\footnote{Two periodic geodesics are called {\it geometrically distinct} if one is not a rotation, inverse or an iterate of the other, i.e., have different images.}. In dimensions up to three, the fundamental group has a strong hold on the geometry of the space, and any such {\it topological} manifold has a unique smooth structure. Since any closed, oriented three manifold is the boundary of a four manifold, one may use this as a bridge to navigate between these dimensions whenever applicable. Techniques of surgery theory has been used succesfully to obtain information in dimension five or higher. \\[0.2cm]
\hf\hf In dimension four the theory is simultaneously rich and mysterious enough for topologists and geometers alike. The vast literature on this subject includes the celebrated results of Freedman and Donaldson, and Taubes' proof of the existence of exotic smooth structures on $\R^4$. It is known that any finitely presented group is the fundamental group of a closed four manifold. The undecidability of the word problem for groups is especially relevant while trying to classify four manifolds. From the early works of Whitehead and Milnor, it is known that the intersection form determines the homotopy type of simply connected, closed four manifold. This intersection pairing can also be interpreted as specific information on Whitehead products of elements of the second homotopy group. \\[0.2cm]
\hf\hf In this paper we reveal a connection between the topology and geometry of closed four manifolds (with finite fundamental group) by answering two questions : {\it how fast do the homotopy groups grow? how many periodic geodesics are there?} The first question is topological in nature and our analysis works for four manifolds which are not necessarily {\it smoothable}. The second question is differential geometric in flavour and is meaningful in the context of Riemannian manifolds. Our proposed solutions to both these queries live in the realm of algebraic topology and homotopy theory. Our main construction is to assign a simply connected, closed, spin 5-manifold to every simply connected, closed, {\it topological} 4-manifold through the construction of certain circle bundles. A crucial ingredient is the knowledge \cite{Sma62} of the diffeomorphism type of such 5-manifolds. The first imprint of this basic idea results in (cf. Theorem \ref{thm45}, Theorem \ref{thm45top}, Theorem \ref{sthomtop}) the following :
\begin{thma}
{\it Let $M$ be a simply connected, closed $4$-manifold with second Betti number $k$.\\
\textup{(i)} The homotopy groups of $M$ are determined by the second Betti number. In particular, if $k\geq 1$ then for $j\geq 3$
\bgd
\pi_j(M)=\pi_j(\#^{k-1} (S^2\times S^3)).
\edd
\textup{(ii)} If $k\geq 1$ then the stable homotopy groups are given by
\bgd
\pi_n^s(M) \cong  (\pi_{n-2}^s)^k \oplus (\pi_{n-3}^s)^{k-1} \oplus \pi_{n-5}^s,
\edd
where $\pi_n^s$ is the $n^{\text{th}}$ stable homotopy groups of spheres.}
\end{thma}
\hf\hf Concerning the topological question recall that a simply connected, closed manifold $M$ is called {\it rationally elliptic} if the sum of the ranks of the rational homotopy groups is finite. A manifold is called {\it rationally hyperbolic} otherwise. Being rationally elliptic imposes constraints strong enough to consider rationally hyperbolic as a generic condition on manifolds. Our second result (cf. Theorem \ref{rnkpi4}, Corollary \ref{MM}) is the following :
\begin{thmb}
{\it Let $M_k$ be a simply connected, closed $4$-manifold with second Betti number $k\geq 3$. For each $j\geq 2$ there exists a polynomial $P_j(x)$ of degree $j$ and top coefficient $1/j$ such that 
\bgd
\textup{rank}\,\pi_{j+1}(M_k)=P_j(k-1).
\edd
As a consequence, $M$ is rationally hyperbolic.}
\end{thmb}
That any $4$-manifold, as above, is rationally hyperbolic is well known (cf. Example 3.8 in \cite{RHT208}) and follows rather tamely from our result. The significance of the theorem lies in its quantitative nature; we provide explicit formula for the ranks of the homotopy groups. As far as the literature of four manifolds is concerned, this is a new result. It is worthwhile pointing out that our scheme of computing these groups is quite general, combining carefully culled inputs from rational homotopy theory and Koszul duality of associative algebras. In fact, we compute the homotopy groups of the associated 5-manifold and this method applies to a larger class of manifolds. \\[0.2cm]
\hf\hf Regarding the geometric question, the study of periodic geodesics has received much attention over the years. It can be traced back to J. Hadamard, H. Poincar\'{e} and G. Birkhoff among others. The natural question, given a compact Riemannian manifold of dimension greater than one, is :
\bgc
{\it Do there exist infinitely many geometrically distinct periodic geodesics?}
\edc
Since periodic geodesics are the critical points of the energy function on the free loop space, Morse theory has been applied in uncovering relationships between periodic geodesics in $M$ and the topology of the free loop space $LM$. The most important result, perhaps, is the theorem of Gromoll and Meyer \cite{GrM69} which states that if the Betti numbers of $LM$ are unbounded then $M$ admits {\it infinitely many} geometrically distinct periodic geodesics in {\it any} metric. Inspired by this result and combining the power of rational homotopy theory, Vigu\'{e}-Poirier and Sullivan \cite{VPSu76} found an affirmative answer to the original question when $M$ has a finite fundamental group and the real cohomology ring of $M$ or any of its covers require at least {\it two} generators. \\[0.2cm]
\hf\hf However, the question for spaces with monogenic cohomology ring is hard even in the classical case of spheres. In the early nineties, using the work of Franks \cite{Fra92}, Bangert \cite{Ban93} (also cf. \cite{Hin93}) proved that the answer is affirmative for {\it any} Riemannian metric on $S^2$. Of particular interest is the result \cite{DuLo10} of Duan and Long that for {\it any} closed, simply connected Riemannian 4-manifold there exist at least {\it two} geometrically distinct periodic geodesics. However, almost all known results either require some restrictions on the curvature or hold only for generic metrics. \\[0.2cm]
\hf\hf A metric is called {\it bumpy} if for each critical point $c$ of the energy function, the circle $S^1(c)$ obtained by rotating $c$, is a nondegenerate critical submanifold of $LM$. Due to Abraham \cite{Abr70} (also cf. \cite{Ano82}) the set of bumpy metrics is {\it dense} in the Frechet space of all the possible metrics on $M$ and can be thought of as a generic condition. We may ask the following question :
\bgc
{\it Do there exist infinitely many geometrically distinct periodic geodesics for a generic metric?}
\edc
An attempt \cite{Kli76} has been made by Klingenberg to settle the question, both for {\it any} metric and for a {\it generic} metric. This has been answered for manifolds with non-zero but finite fundamental group in \cite{BTZ81} while Rademacher \cite{Rad94} proved it for simply connected, compact manifolds. In fact, several other cases are also covered in \cite{BTZ81}.\\[0.2cm]
\hf\hf Let $n_T$ be the number of geometrically distinct periodic geodesics of length at most $T$. In \cite{Gro78} Gromov gave a quantitative version of Gromoll-Meyer's result by giving an estimate for $n_T$ in terms of the Betti numbers $b_i(LM)$'s of the free loop space if the metric is bumpy. Later Ballmann and Ziller \cite{BaZi82} improved this result : for any bumpy metric on a simply connected manifold $M$, there exists {\it positive} constants $a$ and $b$ such that for sufficiently large $T$
\bge\label{BZ}
n_T\geq a \max _{p\leq bT}\, b_p(LM).
\ede
In particular, if the Betti numbers grow exponentially then so does the function $n_T$. Assuming genericity on the part of manifolds, Gromov (cf. \S 5.1 of \cite{RHT208}) asked the following :
\begin{conj}{\bf (Gromov)}\\
{\it Let $M$ be a closed, rationally hyperbolic manifold. For a generic metric on $M$ the function $n_T$ has exponential growth.}
\end{conj}
Our result (cf. Theorem \ref{expLM}, Corollary \ref{expLM4}, Corollary \ref{ReebLM}) is the following :
\begin{thmc}
{\it Let $M$ be a simply connected, closed, smooth $4$-manifold with second Betti number at least $3$. \\
\hf\hf \textup{(1)} For a generic metric on $M$ the number of geometrically distinct periodic geodesics of length at most $\ell$ grows at least exponentially as a function of $\ell$.\\
\hf\hf \textup{(2)} For a generic fibrewise star-shaped hypersurface $\Sigma$ in $T^\ast M$ the number of closed Reeb orbits of length at most $\ell$ grows at least exponentially as a function of $\ell$. }
\end{thmc}
As a consequence of the first part of Theorem C we have a verification of Gromov's conjecture in dimension four. It should be kept in mind that the growth of rational homotopy groups of a space readily controls the growth of the homology of the {\it based} loop space but this is not the case, in general, for the growth of the homology of the {\it free} loop space. This explains why our result is new and also puts it into context. For an explanation of the terminology used in the second part of Theorem C we refer the reader to \S \ref{Reeb}. The result follows without any work once we have Theorem \ref{expLM} and combine this with results from \cite{Hei11}, \cite{MMPa12}.\\[0.5cm]
\hf\hf {\bf Outline of the paper} : In \S \ref{5to4} we explain how to assign a smooth, simply connected, closed, spin $5$-manifold $E$ to any simply connected, closed $4$-manifold $M$ resulting in (i) of Theorem A. We do this first for {\it smooth} $4$-manifolds and the role of the second Stiefel-Whitney class is firmly established. Subsequently we prove this for {\it topological} $4$-manifolds. In \S \ref{sthomgp} we establish (ii) of Theorem A, i.e., the stable homotopy groups of a $M$ are determined by its second Betti number $k$. In \S \ref{geo} we use a result of Lambrechts \cite{Lam01}, which implies that the Betti numbers of $LE$ have exponential growth, as a starting point to prove the first part of Theorem C, i.e., the Betti numbers of $LM$ have exponential growth too. This follows from a delicate and careful analysis of spectral sequences. The subsection \S \ref{ffgp} contains a rather short observation that this exponential growth holds true for $4$-manifolds with finite fundamental group. In \S \ref{Reeb} we use the existing results \cite{MMPa12}, \cite{Hei11} connecting the counting problem for Reeb flows and counting closed geodesics to prove the second part of Theorem C. In \S \ref{growhg} we prove Theorem B, i.e., compute the ranks of the homotopy groups of $M$. We do this in \S \ref{rhg} by computing the homotopy groups of $E$, which is a circle bundle over $M$. Our methods work for a larger class of spaces which includes $E$. Our scheme is to use Milnor-Moore's theorem \cite{MM65} as a common point to merge relevant ideas arising from rational homotopy theory and Koszul duality of associative algebras. Rational hyperbolicity follows trivially from stronger deductions we make on the rational homotopy groups of $M$. Moreover, an explicit formula for the ranks of these groups is given. In \S \ref{details} we prove the details of a key proposition used to prove Theorem B. \\[0.5cm]
\hf\hf {\bf Acknowledgements} : The authors would like to thank Professor John Morgan for some clarifying insights that were useful during the early stages of work on \S \ref{5to4}. The authors are grateful to Professor Felix Schlenk for useful remarks and his insights (cf. \S \ref{Reeb}) on how this work can be used to deduce results on growth for closed Reeb orbits on spherizations.\\[0.5cm]


\section{Associating five manifolds to four manifolds}\label{5to4}

\hf\hf We describe a construction which assigns a smooth, simply connected, closed, spin manifold of dimension five to a closed, simply connected $4$-manifold. The power of the construction arises primarily from :\\
\hf (i) one can identify the diffeomorphism type of this five dimensional manifold;\\
\hf (ii) the $4$-manifold is topological.\\
We start by proving the following theorem for a simply connected, closed smooth $4$-manifold.
\begin{thm}\label{thm45}
For a simply connected, closed smooth $4$-manifold $M$ the homotopy groups are determined by its second Betti number. In particular, if $M$ satisfies $H_2(M;\Z)=\Z^k$ then 
\bgd
\pi_i(M)=\pi_i(\#^{k-1} S^2\times S^3),\,i\geq 3.
\edd 
\end{thm}
Our proof makes crucial use of a classical result of S. Smale. Before embarking on a proof, let's recall Smale's classification result \cite{Sma62}.
\begin{thm}{\bf (Smale '62)}\hf Diffeomorphism classes of smooth, simply connected, closed, spin $5$-manifolds are in bijection with finitely generated abelian groups. The correspondence is given by 
\bgd
N\longrightarrow H_2^\textup{free}(N;\Z)\oplus T
\edd
where $\textup{Tor}(H_2(N;\Z))=T\oplus T$ is direct sum decomposition of the torsion subgroup. Moreover, in the case of simply connected manifolds, any $N$ with $H_2(N;\Z)=\Z^k$ is diffeomorphic to $\#^k S^2\times S^3$.
\end{thm}
We will also need the following.
\begin{lmm}\label{5spin}
Let $M$ be a simply connected, closed, smooth $4$-manifold with $H_2(M;\Z)=\Z^k$. Then there is a fibration $S^1\into E\to M$ such that $E$ is simply connected and spin with $H_2(E;\Z)=\Z^{k-1}$.
\end{lmm}
{\bf Proof.}\hf To construct $E$ we need to deal with two cases - either $M$ is spin or it's not. If $M$ is spin, then choose any element $\alpha$ from a basis of $H^2(M;\Z)=[M,\CP^\infty]$. On the other hand, if $M$ is not spin, then the second Stiefel-Whitney class $w_2$ is the reduction of some integral class $\alpha$, which can be chosen such that it extends to a basis of $H^2(M;\Z)$. In either case, since $\alpha$ corresponds to a homotopy class $[f]$ of maps, let $f:M\to\CP^\infty$ be a representative. In fact, one can assume that $f$ is cellular and $f:M\to\CP^2$. In addition one can assume $f$ to be smooth. Let $E:=f^\ast S^5$ denote the pullback of the Hopf circle bundle over $\CP^2$. Note that $E$ is the associated circle bundle of the line bundle $f^\ast \gamma$ where $\gamma$ is the canonical line bundle over $\CP^2$. By definition, a generator of $H^2(\CP^2;\Z)$ maps to $\alpha$ under $f^\ast$. This implies $c_1(f^*\gamma)=\alpha$. Therefore, one has the Gysin sequence
\bgd
\cdots\longrightarrow H^i(M;\Z)\stackrel{\cup\alpha}{\longrightarrow}H^{i+2}(M;\Z)\stackrel{\pi^\ast}{\longrightarrow}H^{i+2}(E;\Z)\stackrel{\int_{S^1}}{\longrightarrow}H^{i+1}(M;\Z)\longrightarrow\cdots.
\edd  
The Gysin sequence yields exact sequences 
\begin{eqnarray}
\label{spinsm} & 0\longrightarrow H^1(E;\Z) \longrightarrow H^0(M;\Z)\stackrel{\cup\alpha}{\longrightarrow}H^2(M;\Z)\stackrel{\pi^\ast}{\longrightarrow}H^2(E;\Z)\longrightarrow 0 &\\
\label{spinsm2}& 0\longrightarrow H^3(E;\Z) \longrightarrow H^2(M;\Z)\stackrel{\cup\alpha}{\longrightarrow}H^4(M;\Z)\stackrel{\pi^\ast}{\longrightarrow}H^4(E;\Z)\longrightarrow 0 & .
\end{eqnarray}
The first sequence \eqref{spinsm} implies $H^1(E;\Z)=0$ and $H^2(E;\Z)=\Z^{k-1}$. It follows from the universal coefficient theorem that $H_1(E;\Z) = 0$. The long exact sequence of homotopy groups imply
\bgd
0\longrightarrow\pi_2(E)\longrightarrow\pi_2(X)\longrightarrow\Z\longrightarrow\pi_1(E)\longrightarrow 0.
\edd
In particular, $\pi_1(E)$ is a cyclic group and hence abelian. Therefore, $\pi_1(E) \cong H_1(E;\Z)=0$, $E$ is simply connected and hence orientable. Finally, by Poincar\'{e} duality, $H^4(E;\Z)=0$. The equations \eqref{spinsm},\eqref{spinsm2} demonstrate that $H^2(E;\Z)\cong \Z^{k-1} \cong H^3(E;\Z)$. \\
\hf\hf Observe that $E$ is foliated by circles, which carry a nowhere vanishing vector field. This results in a line bundle $L$ over $E$; any such line bundle is trivial as $E$ is simply connected. Since $E$ is the total space of fibration, 
\bgd
TE=\pi^\ast(TM)\oplus L\cong \pi^\ast(TM)\oplus \R.
\edd
Therefore, 
\bgd
w_2(E)=w_2(\pi^\ast(TM))=\pi^\ast(w_2(M))=\pi^\ast(\alpha)=0 
\edd
by construction, and $E$ is spin.  $\hfill\square$\\[0.2cm]
{\bf Proof of Theorem \ref{thm45}.}\hf Let $M_1$ and $M_2$ be two smooth, simply connected and closed $4$-manifolds with the same second Betti number $k$. By Lemma \ref{5spin}, we have two simply connected, closed, spin $5$-manifolds $E_1$ and $E_2$, which are also circle bundles over $M_1$ and $M_2$ respectively. By Smale's result $E_1$ and $E_2$ are mutually diffeomorphic as they are both diffeomorphic to $\#^{k-1} S^2\times S^3$. Therefore, it follows from the long exact sequence of the circle fibrations that
\bgd
\pi_i(M_1)=\pi_i(E_1)=\pi_i(E_2)=\pi_i(M_2),\,i\geq 3.
\edd
Since $\pi_2(M_1)=\pi_2(M_2)=\Z^k$ by Hurewicz isomorphism, we conclude the result. $\hfill\square$ \\[0.2cm]
\hf\hf We want to extend the above construction to {\it topological} $4$-manifolds $X$ which are not necessarily smooth. The argument above fails in such a case as it uses the notion of tangent bundles. In this case define $v_i(X)\in H^i(X;\Z/2)$ to be the cohomology class satisfying 
\bgd
x\cup v_i = \textup{Sq}^i(x)
\edd
for any $x\in H^{4-i}(X;\Z_2)$ \footnote{Note that $v_i$ exists because of Poincar\'{e} Duality with $\Z_2$-coefficients as the cohomology class defined by the linear map $x\mapsto \langle Sq^i(x),[X]\rangle$.}. We then define $w_i(X)\in H^i(X;\Z_2)$ by the equation\footnote{The reader should recall Wu's formula which presents the total Stiefel-Whitney class as the Steenrod square of the Wu class.}
\bgd
1+w_1+w_2+w_3+w_4= (1+\textup{Sq}^1+\textup{Sq}^2+\textup{Sq}^3)(1+v_1+v_2+v_3+v_4).
\edd
If $X$ is simply connected we have $v_1=0$, whence $w_1=0$. Therefore, the above equation implies that $w_2=v_2$. We prove an analogue of Lemma \ref{5spin}.
\begin{lmm}\label{5spintop}
Let $X$ be a simply connected, closed, topological $4$-manifold with $H_2(X)\cong \Z^k$. Then there exists a circle bundle $S^1\hookrightarrow E \rightarrow X$ such that $E$ is a smooth $5$-manifold which is simply connected and spin. In addition, $H_2(E)\cong \Z^{k-1}$.
\end{lmm}  
{\bf Proof.}\hf We have two cases according to $w_2(X)=0$ and $w_2(X)\neq 0$. If $w_2(X)=0$, then choose any primitive element $\alpha$ of $H^2(X;\Z)=[X,\CP^\infty]\cong \Z^k$, i.e., $\alpha$ is not a multiple of any other element. On the other hand, if $w_2(X)\neq 0$ then $w_2(X)$ is of the form $\alpha~(\mbox{mod}~2)$ for some primitive element $\alpha$ of $H^2(X;\Z)$. In either case,  $\alpha$ corresponds to a homotopy class $[f]$ for $f:X\to\CP^\infty$ which may assumed to be cellular, i.e., $f:X\to\CP^2$. Let $E:=f^\ast S^5$ denote the pullback of the circle bundle over $\CP^2$. It follows from the first part of the argument in the proof of Lemma \ref{5spin} that $E$ is closed, simply connected $5$-manifold with $H^2(E;\Z)=\Z^{k-1}$.\\
\hf We know that \cite{KiSi77} a topological manifold $Y$ admits a PL-structure if and only if the Kirby-Siebenmann class $e(E)\in H^4(Y;\Z_2)$ vanish. In our case, $H^4(E;\Z_2)=0$ and $E$ admits a PL-structure. In follows from the work of various people (cf. \cite{Mil11} for a detailed history) that in dimension $7$ or lower every PL-structure possesses a compatible smooth structure. This smooth structure is unique up to diffeomorphism for dimension $6$ or lower. In particular, $E$ is a smooth manifold.\\
\hf It remains to verify that $E$ is spin. Let $T(f^\ast \gamma)$ denote the Thom space of the bundle $f^\ast \gamma$. The Thom isomorphism theorem states that $H^{2+k}(T(f^\ast \gamma))\cong H^k(X)$ and the isomorphism is induced by multiplication with the Thom class of the bundle. The Thom space can be identified with $\mathcal{D}(f^\ast \gamma)/S(\gamma)$ where $\mathcal{D}(f^\ast \gamma)$ is the disk bundle associated to $f^\ast \gamma$ which is homotopy equivalent to $X$, and $S(\gamma)=E$. In cohomology with $\Z_2$-coefficients we have the long exact sequence
\bgd
\cdots\longrightarrow \widetilde{H}^i(T(f^\ast \gamma);\Z_2)\longrightarrow H^i(X;\Z_2)\stackrel{\pi^\ast}{\longrightarrow}H^i(E;\Z_2)\longrightarrow \widetilde{H}^{i+1}(T(f^\ast \gamma);\Z_2)\longrightarrow\cdots.
\edd  
\hf We will prove that $\textup{Sq}^2(a)=0$ for all $a\in H^*(E;\Z_2)$. Note that 
\bgd
H^k(E;\Z_2)=\left \{\begin{array}{rl} 
                \Z_2 &\mbox{if}~k=0,5 \\
                \Z_2^{k-1} &\mbox{if}~k=2,3  \\
                0    &\mbox{otherwise}. 
\end{array}\right.
\edd 
Hence, $\textup{Sq}^2:H^i(E;\Z_2)\rightarrow H^{i+2}(E;\Z_2)$ maps $a\in H^i(E;\Z_2)$ to $0$ for all $i\neq 3$. Note that from the Thom isomorphism 
\bgd
H^k(T(f^\ast \gamma);\Z_2)=\left \{\begin{array}{rl} 
                \Z_2 &\mbox{if}~k=0,2,6 \\
                \Z_2^k &\mbox{if}~k=4  \\
                0    &\mbox{otherwise}                   
\end{array}\right.
\edd
The long exact sequence above implies the map $H^j(E;\Z_2) \rightarrow H^{j+1}(T(f^\ast \gamma);\Z_2)$ is injective for $j=3,5$. Since this comes from the boundary map of the long exact sequence it preserves the $\textup{Sq}^2$-action, i.e., we have the commuting diagram
\bgd
\xymatrix{
H^3(E;\Z_2)\ar[r]^-{i}\ar[d]_-{\textup{Sq}^2} & H^4(T(f^\ast \gamma);\Z_2)\ar[d]^-{\textup{Sq}^2}\\
H^j(E;\Z_2)\ar[r]^-{i} & H^{4}(T(f^\ast \gamma);\Z_2)}
\edd
Therefore, it suffices to show that $\textup{Sq}^2$ is zero on the image $i(H^3(E;\Z_2))$. Let $\phi\in H^2(T(f^\ast \gamma);\Z_2)$ be the Thom class of $f^\ast \gamma$; this maps to $\alpha$ (mod $2$) under the map $ H^2(T(f^\ast \gamma);\Z_2) \rightarrow H^2(X;\Z_2)$. We call this class $\alpha$ as well.\\
\hf Suppose first that $w_2(X)=0$, which means $\textup{Sq}^2(\lambda)=0$ for all $\lambda \in H^2(X;\Z_2)$. The group $H^4(T(f^\ast \gamma);\Z_2)$ consists of products $\lambda \phi$ for $\lambda \in H^2(X;\Z_2)$. Then 
\bgd
\textup{Sq}^2(\lambda \phi)= \textup{Sq}^2(\lambda)\phi + \textup{Sq}^1(\lambda)\textup{Sq}^1(\phi)+\lambda \textup{Sq}^2(\phi)= (\lambda \alpha) \phi.
\edd
Thus, $\textup{Sq}^2(\lambda \phi)\neq 0$ if and only if $\lambda \alpha \neq 0$. The latter condition implies that under the map
\bgd
H^4(T(f^\ast L);\Z_2) \longrightarrow H^4(X;\Z_2)
\edd
the element $\lambda \phi$ is sent to $\lambda \alpha$. Consequently, $\lambda \phi$ does not lie in the image of $H^3(E;\Z_2)$, whence $\textup{Sq}^2$ vanishes on $i(H^3(E;\Z_2)$. \\
\hf Now suppose $w_2(X)\neq 0$. Then $\alpha^2=\textup{Sq}^2(\alpha)\neq 0$ since multiplication by $\alpha$ is the operation $\textup{Sq}^2$ which by assumption is non-zero. In this case
\bgd
\textup{Sq}^2(\lambda \phi)= \textup{Sq}^2(\lambda) \phi + \textup{Sq}^1(\lambda)\textup{Sq}^1(\phi) + \lambda \textup{Sq}^2(\phi) = \lambda \alpha  \phi +\lambda \alpha \phi = 0.
\edd
\hf We conclude that the action of $\textup{Sq}^2$ on $H^*(E;\Z_2)$ is $0$. This implies that $v_2(E)=0$. We already have that $v_1(E)=0$ since $E$ is simply connected. Therefore, $w_2(E)=0$. By Wu's formula, we have $w_2(E)=w_2(TE)$, and therefore, $E$ is spin. $\hfill\square$\\[0.2cm]
We can now conclude an analogue of Theorem \ref{thm45}.
\begin{thm}\label{thm45top}
For a simply connected, closed $4$-manifold $X$ the homotopy groups are determined by its second Betti number. In particular, if $X$ satisfies $H_2(X;\Z)=\Z^k$ then 
\bgd
\pi_i(X)=\pi_i(\#^{k-1} S^2\times S^3),\,i\geq 3.
\edd 
\end{thm}
{\bf Proof.}\hf Let $X_1$ and $X_2$ be two simply connected, closed $4$-manifolds with the same second Betti number $k$. By Lemma \ref{5spintop}, we have two spin (i.e, smooth, simply connected, closed and spin) $5$-manifolds $E_1$ and $E_2$, which are also circle bundles over $X_1$ and $X_2$ with isomorphic homology groups. By Smale's result $E_1$ and $E_2$ are mutually diffeomorphic and they are both diffeomorphic to $\#^{k-1} S^2\times S^3$. Therefore, we conclude the result as before.$\hfill\square$
\begin{eg}{\bf (Four manifolds of rank two)}\hf Let's see explicitly what happens for $S^2\times S^2$ and $\CP^2\#\overline{\CP}^2$. Recall that the principal $S^1$-bundle $S^5\stackrel{\pi}{\longrightarrow}\CP^2$ corresponds to the first Chern class $\mathpzc{c}\in H^2(\CP^2;\Z)$ of the tautological line bundle. Equivalently, if we want to kill the obstruction class to $\CP^2$ being spin, then we construct a principal circle bundle via the correspondence between $H^2(X;\Z)$ and $[X,\CP^\infty]$. Moreover, $w_2(\CP^2)=w_2(\overline{\CP}^2)$ since we're working over $\Z_2$. Therefore, we have a corresponding bundle $S^1\into S^5\to\overline{\CP}^2$ which kills $w_2(\overline{\CP}^2)$ when pulled back to $S^5$. Let $D$ denote a disk of radius $1$ around $[1:0:0]\in\CP^2$, viz.,
\bgd
D:=\left\{[1:z_1:z_2]\,|\,|z_1|^2+|z_2|^2\leq 1\right\}.
\edd
This implies that
\begin{eqnarray*}
\pi^{-1}(D) & := & \left\{(\lambda\cos\theta,\lambda\sin\theta, x_3,\ldots, x_6)\,|\,\theta\in[0,2\pi),\,x_3^2+x_4^2+x_5^2+x_6^2=1-\lambda^2\leq\lambda^2\neq 0\right\}\\
 & = & \left\{(x_1,x_2,x_3\ldots,x_6)\,|\,x_3^2+x_4^2+x_5^2+x_6^2\leq{\textstyle \frac{1}{2}},\,x_1^2+\ldots+ x_6^2=1\right\} \cong S^1\times D_4
\end{eqnarray*}
while the complement of $\pi^{-1}(D)$ is 
\bgd
\{(x_1,\ldots, x_6)\,|\,x_1^2+x_2^2\leq {\textstyle \frac{1}{2}},\,x_1^2+\ldots+ x_6^2=1\}\cong D_2\times S^3.
\edd
Therefore, $S^5$ can be written as the union of $S^1\times D_4$ and $D_2\times S^3$. \\
\hf\hf Let $X:=S^5\setminus \pi^{-1}(D)$. Then we have fibrations
\bgd
S^1\into X\to \CP^2\setminus D,\,\,\,\,\,\,\,S^1\into X\to \overline{\CP}^2\setminus D.
\edd
Let $\widetilde{X}$ be the space obtained by gluing two copies of $X$ together by the identity along its boundary. It's clear that $\widetilde{X}=S^2\times S^3$ and the two fibrations above blend into a fibration
\bgd
S^1\into S^2\times S^3\to \CP^2\#\overline{\CP}^2.
\edd
The total space $S^2\times S^3$ is also a circle bundle over $S^2\times S^2$. However, the Whitehead product on the homotopy groups help distinguish $S^2\times S^2$ from $\CP^2\#\overline{\CP}^2$. 
\end{eg}


\subsection{Stable homotopy groups}\label{sthomgp}

\hf\hf We exploit the circle fibre bundle introduced in \S \ref{5to4} to gather information about stable homotopy groups of a simply connected, closed $4$-manifold $M$ in terms of the second Betti number $k$. Let $\pi_n^s$ denote the $n^\text{th}$ stable homotopy groups of spheres. We can directly compute the stable homotopy group $\pi_n^s(M)$ of $M$ as follows.  
\begin{thm}\label{sthomtop}
Let $M$ be a simply connected, closed $4$-manifold with second Betti number $k\geq 1$. Then
\bge\label{stpi}
\pi_n^s(X) \cong  (\pi_{n-2}^s)^k \oplus (\pi_{n-3}^s)^{k-1} \oplus \pi_{n-5}^s.
\ede
\end{thm}
{\bf Proof.}\hf We use the fibration $S^1\rightarrow \#^{k-1}(S^2\times S^3) \rightarrow M$. On passing to suspension spectra in the stable homotopy category, we obtain a fibre sequence 
\bgd
\Sigma^\infty S^1 \rightarrow \Sigma^\infty \#^{k-1}(S^2\times S^3) \rightarrow \Sigma^\infty M.
\edd
In the stable homotopy category, fibre sequences are also cofibre sequences which yields a long exact sequence of homotopy groups. 
\bgd
\cdots\rightarrow \pi_n^s(\#^{k-1}(S^2\times S^3))\rightarrow \pi_n^s(M) \rightarrow \pi_{n-1}^s(S^1) \rightarrow \pi_{n-1}^s(\#^{k-1}(S^2\times S^3))\rightarrow\cdots .
\edd
\hf The first observation is that $\#^{k-1}(S^2\times S^3)$ has a CW-structure with trivial $1$-skeleton, and hence the map $S^1 \to \#^{k-1}(S^2\times S^3)$ is null-homotopic. Therefore, the the map induced on spectra is trivial. Since fibre sequences are cofibre sequences in the stable homotopy category, $\Sigma^\infty M$ being the cofibre of the map $\Sigma^\infty S^1\to\Sigma^\infty \#^{k-1}(S^2\times S^3)$ becomes the wedge of $\Sigma^\infty \#^{k-1}(S^2\times S^3)$ and $S^2$. Therefore,
\bgd
\pi_n^s(M) \cong \pi_n^s(\#^{k-1}(S^2\times S^3)) \oplus \pi_{n-1}^s(S^1) \cong \pi_n^s(\#^{k-1}(S^2\times S^3)) \oplus \pi_{n-2}^s.
\edd
We know that suspension of the Whitehead product is trivial. The space $\#^{k-1}(S^2\times S^3)$ has a CW-complex structure with $4$-skeleton a wedge of $(k-1)$ two spheres ($a_i$'s) and $(k-1)$ three spheres ($b_j$'s), and a single $5$-cell attached along the map $[a_1,b_1] +\ldots +[a_{k-1},b_{k-1}]$. After a suspension we have 
\bgd
\Sigma \#^{k-1}(S^2\times S^3) \simeq \vee^{k-1} (\Sigma S^2 \vee \Sigma S^3) \vee \Sigma S^5 \simeq  \vee^{k-1}(S^3\vee S^4) \vee S^6.
\edd
This implies 
\bgd
\pi_n^s(\#^{k-1}(S^2\times S^3)) \cong \pi_n^s(S^2)^{k-1} \oplus \pi_n^s (S^3)^{k-1} \oplus \pi_n^s(S^5) \cong (\pi_{n-2}^s)^{k-1} \oplus (\pi_{n-3}^s)^{k-1} \oplus \pi_{n-5}^s.
\edd
Therefore, 
\bgd
\pi_n^s(M) \cong (\pi_{n-2}^s)^{k-1} \oplus (\pi_{n-3}^s)^{k-1} \oplus \pi_{n-5}^s \oplus \pi_{n-2}^s\cong  (\pi_{n-2}^s)^k \oplus (\pi_{n-3}^s)^{k-1} \oplus \pi_{n-5}^s 
\edd
as claimed. $\hfill\square$ \\[0.2cm]
\hf\hf We can extend the above argument to the case where $M$ is closed and has {\it finite} fundamental group. The second homotopy group is $\pi_2(M)=\Z^k$ for some $k$. This follows from the fact that $\widetilde{M}$ is a simply connected, closed $4$-manifold with $\pi_2(M)=\pi_2(\widetilde{M})=\Z^k$. 
\begin{cor}
\label{sthomtop2}
Let $M$ be a closed $4$-manifold with a fundamental group of order $m$ and $\pi_2(M)=\Z^k$ with $k\geq 1$. Then
\bgd
\pi_n^s(M) \cong  (\pi_{n-2}^s)^k \oplus (\pi_{n-3}^s)^{k-1} \oplus \pi_{n-5}^s\oplus (\pi_{n-1}^s)^{m-1}.
\edd
\end{cor}
{\bf Proof.}\hf The covering map $\pi:\widetilde{M}\to M$ is $m$ to $1$. In the stable homotopy category we have a (fibre and cofibre) sequence 
\bgd
\cdots\rightarrow \pi_n^s(\widetilde{M})\rightarrow \pi_n^s(M) \rightarrow \pi_{n-1}^s(\{m\,\textup{points}\}) \rightarrow \pi_{n-1}^s(\widetilde{M})\rightarrow\cdots .
\edd
The inclusion map of a fibre consisting of $m$ points into $\widetilde{M}$ is nullhomotopic. Therefore,
\bgd
\pi_n^s(M) \cong \pi_n^s(\widetilde{M}) \oplus \pi_{n-1}^s(\{m\,\textup{points}\}).
\edd
Since $\Sigma^j(\{m\,\textup{points}\})=\vee^{m-1}S^j$, combining this with \eqref{stpi} we have
\bgd
\pi_n^s(M) \cong \pi_n^s(\widetilde{M}) \oplus (\pi_{n-1}^s)^{m-1}
\edd
which finishes the proof. $\hfill\square$


\section{How many periodic geodesics are there?}\label{geo}

\hf\hf The question of how many geometrically distinct periodic geodesics exist on a simply connected closed, Riemannian manifold $M$ is classical. It follows from results of Gromoll-Meyer \cite{GrM69} that the answer is linked to the Betti numbers of the {\it free loop space} $LM:=C^0(S^1,M)$. We should point out that the Morse theory techniques are usually applied to the space $\mathcal{L}M$ of {\it closed loops of class $H^1$}. It is known that $LM$ is homotopy equivalent to its geometric cousin $\mathcal{L}M$ which is a Hilbert manifold. One studies the critical points of the energy function on $\mathcal{L}M$. Periodic geodesics are precisely the critical points but these are not necessarily geometrically distinct. Consequently, a careful analysis of the multiplicity needs to be done. 
\begin{thm}{\bf (Gromoll - Meyer '69)}\hf Let $M$ be a simply connected, closed Riemannian manifold. If there exists a field $\mathbf{k}$ such that the sequence of ranks of the homology groups of $LM$ are unbounded, then for each Riemannian metric there exists infinitely many geometrically distinct periodic geodesics on $M$.
\end{thm}
Towards an answer to the question of when manifolds have unbounded growth in the homology of its loop space, V-Poirier and Sullivan \cite{VPSu76} had the following to say :
\begin{thm}{\bf (V-Poirier - Sullivan '76)}\hf If $M$ is a closed Riemannian manifold with a finite fundamental group such that $M$ or any of its covers require at least two generators for its cohomology then the ranks of $H_\ast(LM)$ is unbounded. In particular, there exists infinitely many geometrically distinct periodic geodesics on $M$.
\end{thm}
This, for instance, tells us that any smooth, simply connected, closed Riemannian $4$-manifold $M^4$ with $\chi(M^4)\geq 4$ admits infinitely many geometrically distinct periodic geodesics. \\
\hf\hf To proceed further, let us recall an important but elementary notion.
\begin{defn}\label{expseq}
Let $\{b_n\}_{n\geq 1}$ be a sequence of non-negative numbers. We say that this sequence has {\it exponential growth} if there exists $\lambda>0, C>1$ and $N\geq 1$ such that
\bgd
\textup{max}\{b_1,\cdots,b_n \}\geq \lambda C^n\,\,\textup{for all}\,\,n\geq N.
\edd
\end{defn}
\begin{rem}\label{exp2}
The definition of exponential growth is equivalent to an alternate definition where one requires
\bgd
b_1+\cdots+ b_n\geq \lambda C^n
\edd
for $\lambda>0, C>1$ and $n\geq N$ for some $N$. More importantly, the defintion above does nt distinguish between exponential growth and super-exponential growth, for instance. We shall mean at least exponential growth when we say exponential growth.
\end{rem}
We are ready to state the following \cite{Lam01} useful result.
\begin{thm}{\bf (Lambrechts '01)}\hf Let $M_1$ and $M_2$ be two simply connected, closed manifolds of dimension $n$. Suppose that $\mathbf{k}$ is a field such that neither $H^\ast(M_1;\mathbf{k})$ nor $H^\ast(M_2;\mathbf{k})$ is isomorphic to $H^\ast(S^n;\mathbf{k})$. Then\\
\hf \textup{(1)} the sequence $\{b_n\}_{n\geq 1}$ given by
\bgd
b_n:=\dim H_n(L(M_1\# M_2);\mathbf{k})
\edd
is unbounded;\\
\hf \textup{(2)} we have the following dichotomy :\\
\hf\hf \textup{(a)} either both of the cohomology algebras $H^\ast(M_i;\mathbf{k})$ are monogenic and the sequence $\{b_n\}$ has a linear growth;\\
\hf\hf \textup{(b)} or at least one $H^\ast(M_i;\mathbf{k})$ is not a monogenic algebra and the sequence $\{b_n\}$ has an exponential growth.
\end{thm}
As an iterated application of the result above we have the following :
\begin{cor}\label{E5exp}
For $k\geq 0$ let $\mathpzc{E}_k:=\#^k S^2\times S^3$ be the $5$-manifold which is the connected sums of $k$ copies of $S^2\times S^3$. If $k\geq 2$ then the Betti numbers of the free loop space $L\mathpzc{E}_k$ grow exponentially.
\end{cor}
\hf\hf Let $M^4_{k+1}$ be any simply connected, closed, smooth $4$-manifold with $\chi(M^4_{k+1})=k+3$ (or equivalently, $H_2(M^4_{k+1})=\Z^{k+1}$). In light of the connection (cf. \S \ref{5to4}) between $M^4_{k+1}$ and $\mathpzc{E}_k$, a natural question is the following :
\bgc
{\it Does the sequence $\{\dim H_n(LM^4_{k+1};\Q)\}_{n\geq 1}$ grow at least exponentially?}
\edc
For the ease of reading and brevity, we shall fix $k\geq 2$ for now and denote $M^4_{k+1}$ by $M$. We will denote $\mathpzc{E}_k$ simply by $\mathpzc{E}$.\\
\hf\hf We begin with the fibration $LS^1\to L\mathpzc{E}\to LM$. Notice that for any space $X$, $LX$ is itself a fibration over $X$ with fibre $\Omega X:=\textup{Map}_\ast(S^1,X)$, the based loops in $X$. Moreover, there is a section $s:X\to LX$ by sending a point to the constant loop at that point. Therefore,
\bgd
\pi_i(LX)=\pi_i(X)\oplus\pi_i(\Omega X)=\pi_i(X)\oplus \pi_{i+1}(X).
\edd
In particular, this means
\bgd
\pi_1(LM)=\pi_2(M)=H_2(M)=\Z^{k+1}.
\edd
Moreover, $\pi_1(LM)$ acts on the fibre $LS^1\cong S^1\times \Z$ (of the fibration $LS^1\into L\mathpzc{E}\to LM$). Applying the Serre spectral sequence to this fibration, we get $E^2$-page :
\bgc
\begin{tabular}{|ccccc}
$\vdots$ & $\vdots$ & $\vdots$ & $\vdots$ & $\iddots$\\[5pt]
$0$ & $0$ & $0$ & $0$ & $\cdots$ \\[10pt]
$H_0(LM;\mathcal{R}_1)$ & $H_1(LM;\mathcal{R}_1)$ & $H_2(LM;\mathcal{R}_1)$ & $H_3(LM;\mathcal{R}_1)$ & $\cdots$\\[15pt]
$H_0(LM;\mathcal{R}_0)$ & $H_1(LM;\mathcal{R}_0)$ & $H_2(LM;\mathcal{R}_0)$ & $H_3(LM;\mathcal{R}_0)$ & $\cdots$\\[8pt]
\hline
\end{tabular}
\edc
Observe that the coefficients are given by
\bgd
\mathcal{R}_0=H_0(LS^1),\,\mathcal{R}_1=H_1(LS^1).
\edd
A generator $\mathbbm{1}\in\pi_1(LM)$ acts as a homotopy equivalence on $LS^1$, and induces the same action on $H_0(LS^1;\Z)$ and $H_1(LS^1;\Z)$. Both of these are isomorphic to the same abelian group $\Z[\Z]$ but are to be thought of as {\it twisted coefficients} in the spectral sequence above. Moreover, it is evident from the structure of the $E^2$-page that $E^3=E^\infty$. Therefore, it suffices to analyze this page as we're only interested in the additive information.
\begin{lmm}
The homology of $LM$ with twisted coefficients $\mathcal{R}_i,i=0,1$ is isomorphic to the homology of a cover $\mathcal{L}$ of $LM$.
\end{lmm}
{\bf Proof.}\hf The coefficient system $\mathcal{R}_i$ has the structure of a $\Z[\pi_1(LM)]$-module. Let $G:=\pi_1(LM)$ for notational convenience. In general, $H_\ast(LM;\mathcal{R})$ can be computed (or even identified) as the homology of the chain complex
\bgd
C_\ast(\widetilde{LM};\Z)\otimes_{\Z[G]} \mathcal{R}.
\edd
Let us denote the action of $G$ on $H_0(LS^1)$ by $\rho:G\to \textup{Aut}(H_0(LS^1))$. Let $\mathcal{L}$ be the cover of $LM$ that corresponds to $\textup{ker}\,\rho$. Then it follows that
\bgd
C_\ast(\widetilde{LM};\Z)\otimes_{\Z[G]} \mathcal{R}\cong C_\ast(\mathcal{L};\Z)\otimes_{\Z[G/\textup{ker}\,\rho]} \mathcal{R}.
\edd
We observe that 
\bgd
\textup{ker}\,\rho\cong\Z^{k-1},\,\,\,\,(G/\textup{ker}\,\rho)\cong \Z.
\edd
In fact, the action of $G$ on $H_0(LS^1)$ and $H_1(LS^1)$ is identical. Therefore, 
\bgd
H_\ast(LM;\mathcal{R}_i)\cong \textup{homology of}\,\left(C_\ast(\mathcal{L};\Z)\otimes_{\Z[\Z]}\mathcal{R}_i\right)\cong H_\ast(\mathcal{L};\Z)
\edd
as $\mathcal{R}_i\cong \Z[\Z]$ as $\Z[\Z]$-modules. $\hfill\square$\\[0.2cm]
\hf\hf To proceed further, let us denote the maps on the $E^2$-page by
\bgd
d:H_n(\mathcal{L};\Z)\longrightarrow H_{n-2}(\mathcal{L};\Z).
\edd
Then it follows that
\bgd
H_n(L\mathpzc{E};\Z)\cong \frac{H_{n-1}(\mathcal{L};\Z)}{d(H_{n+1}(\mathcal{L};\Z))}\oplus \left(\textup{ker}:H_n(\mathcal{L};\Z)\longrightarrow H_{n-2}(\mathcal{L};\Z)\right).
\edd
Let $\beta_n$ be the rank of $H_n(\mathcal{L};\Z)$ while $b_n$ be the rank of $H_n(L\mathpzc{E};\Z)$. We have the inequality
\bgd
b_n\leq \beta_n+\beta_{n-1}.
\edd
It follows either from Defintion \ref{expseq} or from Remark \ref{exp2} that $\{\beta_n\}_{n\geq 1}$ has exponential growth if $\{b_n\}_{n\geq 1}$ does. Combined with Corollary \ref{E5exp} we conclude the following :
\begin{cor}\label{expL}
The Betti numbers $\{\textup{dim}\,H_n(\mathcal{L};\Z)\}_{n\geq 1}$ has exponential growth.
\end{cor}
\hf\hf In order to arrive at a statement involving the growth of the Betti numbers of $LM$ we observe that $\mathcal{L}\to LM$ is a principal $\Z$-bundle and is classified by a (homotopy class of) maps $\varphi:LM\to B\Z$. But as $\Z\simeq \Omega S^1$ we may also think of $\mathcal{L}$ as the pullback of $\textup{exp}:\R \to S^1$ via some map $\varphi:LM\to S^1$. 
\begin{thm}\label{expLM}
The Betti numbers $\{\textup{dim}\,H_n(LM;\Z)\}_{n\geq 1}$ has exponential growth.
\end{thm}
{\bf Proof.}\hf We turn the map $\varphi:LM\to S^1$ into a fibration where the fibre has the homotopy type of $\mathcal{L}$. We have the Serre spectral sequence for $\mathcal{L}\to LM\to S^1$, where $\pi_1(S^1)=\Z$ acts on the fibre as the usual Deck transformation. The Serre spectral sequence at the $E^2$-page is all we have and we get
\bgd
H_j(LM;\Z)=H_0(S^1;H_j(\mathcal{L};\Z))\oplus H_1(S^1;H_{j-1}(\mathcal{L};\Z)).
\edd
By Proposition \ref{1id}, the generator of $\pi_1(S^1)$ induces the identity map on $H_\ast(\mathcal{L};\Z)$ . The spectral sequence then reads
\bgd
H_j(LM;\Z)=H_j(\mathcal{L};\Z)\oplus H_{j-1}(\mathcal{L};\Z).
\edd
Since the Betti numbers of $\mathcal{L}$ has exponential growth (cf. Corollary \ref{expL}), the same hold for the Betti numbers of $LM$. $\hfill\square$
\begin{cor}\label{expLM4}
For a generic metric on $M$ the number of geometrically distinct periodic geodesics of length at most $\ell$ grows exponentially as a function of $\ell$.
\end{cor}
{\bf Proof.}\hf The result follows from Ballmann-Ziller's estimate \cite{BaZi82} (cf. \eqref{BZ} in \S \ref{intro}) and Theorem \ref{expLM}. $\hfill\square$
\begin{rem}
The growth of the Betti numbers of $L(S^2\times S^2)$ or $L(\CP^1\#\overline{\CP}^2)$ is linear. This is related to the fact that $S^2\times S^3$, the associated $5$-manifold to both these $4$-manifolds, is rationally elliptic. It is known \cite{Hin93} that the function $n_T$ grows at least as fast $\frac{T}{\log T}$ for any metric on $S^2$. The function $n_T(S^2\times S^2)$ grows at least as fast as $(T/\log T)^2$ for any product metric on $S^2\times S^2$. 
\end{rem}
\hf\hf The proof of Theorem \ref{expLM} is contingent on the fact that the local coefficients in the Serre spectral sequence can be replaced by untwisted coefficients. Towards that end we have :
\begin{prpn}\label{1id}
Let $\mathbbm{1}$ denote a generator of $\Z\subset \pi_1(LM)$ which correponds to the covering $\mathcal{L}$. We write 
\bgd
\pi_1(LM)=\Z^{k+1}=\langle \mathbbm{1} \rangle \oplus G,
\edd
where $G=\Z^k$. The action of $\mathbbm{1}$ on $\widetilde{LM}$ induces the identity map at the level of homology. In particular, the induced $G$-equivariant map
\bgd
\mathbbm{1}^G_\ast: H_\ast(\mathcal{L};\Z)\longrightarrow H_\ast(\mathcal{L};\Z)
\edd
is the identity map.
\end{prpn}
{\bf Proof.}\hf Let us suppose that 
\bge\label{bbm1}
\mathbbm{1}:\widetilde{LM}\longrightarrow \widetilde{LM}
\ede
induces the identity map on homology. Since $G$ (the direct summand of $\langle \mathbbm{1}\rangle$ in $\pi_1(LM)$) acts freely on $\widetilde{LM}$ and $\pi_1(LM)$ is abelian, the map \eqref{bbm1} is a $G$-equivariant map. If $\mathbbm{1}_\ast=\textup{Id}$ then $\mathbbm{1}^G_\ast=\textup{Id}$, where $\mathbbm{1}^G_\ast$ is the induced map on the $G$-equivariant homology. To see this claim recall that the equivariant homology is the homology of $H_\ast(EG\times_G X)$. Hence there is a spectral sequence\footnote{The Serre spectral sequence of the fibration $X\rightarrow EG\times_G X \rightarrow BG$ computes the requisite homology.} calculating the this homology whose $E^2$-page is $H_p(G;H_q(X))$ and which converges to $H_{p+q}(EG\times_G X)$. A $G$-equivariant map $f$ induces a self-homomorphism of this spectral sequence which is identity on the $E^2$-page. Therefore, it induces the identity map on the $E^\infty$-page, and hence, on the equivariant homology. \\
\hf Since $\mathbbm{1}^G_\ast$ is the map on the $G$-equivariant homology of $\widetilde{LM}$ and the action is free,
\bgd
H^G_\ast(\widetilde{LM};\Z)\cong H_\ast(\widetilde{LM}/G;\Z)=H_\ast(\mathcal{L};\Z).
\edd
This proves the second claim in the proposition.\\
\hf To prove that $\mathbbm{1}_\ast$ (in \eqref{bbm1}) is the identity on homology of $\widetilde{LM}$, we notice that $LM$ fibres over $M$ with fibre $\Omega M$, the based loop space. This fibration has a section $s$ given by constant loops. We may replace each fibre $\Omega_p M$ (the loops in $M$ based at $p\in M$) by its universal cover $\widetilde{\Omega_p M}$ based at $p$. We may do so in a continuous fashion due to the presence of the section $s:M\to LM$; we construct the path space and then impose the necessary equivalence relations to construct the universal cover. We get a new fibration (with $\widetilde{LM}$ as the total space)
\bge\label{tildeLM}
\widetilde{\Omega M}\into \widetilde{LM}\longrightarrow M.
\ede
Recall that $\mathbbm{1}$ is a Deck transformation of the covering map $\widetilde{LM}\to LM$. Since $LM$ fibres over $M$, the map $\mathbbm{1}$ is one of fibrations, mapping \eqref{tildeLM} to itself and is the identity map on $M$. The map $\mathbbm{1}_\ast$ induces a map of spectral sequences. By Lemma \ref{topmon} applied to $G=\Omega_p M$ with $e$ as the constant loop at $p\in M$, we conclude that $\mathbbm{1}_\ast$ acts as identity on $H_\ast(\widetilde{\Omega M};\Z)$. It is already identity on $H_\ast(M;\Z)$ and a comparison of the map between spectral sequences in homology tells us $\mathbbm{1}_\ast$ is identity on the $E^2$-page. This persists till the $E^\infty$-page and it follows that $\mathbbm{1}_\ast$ is the identity map on the homology of $\widetilde{LM}$. $\hfill\square$\\[0.2cm]
\hf\hf To complete the proof of the proposition above we need the following :
\begin{lmm}\label{topmon}
Let $G$ be a topological monoid with universal cover $\widetilde{G}$. Let $e\in G$ be an element such that $e\cdot e=e$ and $e$ is an identity up to homotopy. If\, $\Z=\langle \mathbbm{1}\rangle\subset \pi_1(G,e)$ then 
\bgd
\mathbbm{1}_\ast:H_\ast(\widetilde{G};\Z)\longrightarrow H_\ast(\widetilde{G};\Z)
\edd
is the identity map.
\end{lmm}
{\bf Proof.}\hf Fix an element $e\in G$. We think of $\widetilde{G}$ as the universal cover arising from the path space associated to $G$, i.e., elements of $\widetilde{G}$ are equivalence classes of continuous paths $\gamma:[0,1]\to G$ such that $\gamma(0)=e$. We say $\gamma$ and $\eta$ are equivalent if $\gamma(1)=\eta(1)$ and $\gamma\ast\eta^{-1}$ is null homotopic. We also denote by $\widetilde{G}_g$ the path space associated to $G$ with $g$ as the chosen point for all paths to start at. In this notation, $\widetilde{G}=\widetilde{G}_e$. Choose a representative $\kappa:(S^1,1)\to (G,e)$ for $\mathbbm{1}$. Then the action of $\mathbbm{1}$ is given by
\bgd
\mathbbm{1}:\widetilde{G}\to \widetilde{G},\,\,\,\,[\gamma]\mapsto [\kappa\ast\gamma].
\edd
Define a continuous family of maps
\bgd
\mathbbm{1}^t:\widetilde{G}\longrightarrow\widetilde{G},\,\,\,\,[\gamma]\mapsto \big[(\kappa|_{[0,t]}\cdot e)\ast(\kappa(t)\cdot\gamma)\big].
\edd
Since $e$ is an identity up to homotopy, we have
\begin{eqnarray*}
\mathbbm{1}^0\big([\gamma]\big)=[(e\cdot e)\ast (e\cdot \gamma)]=[e\ast \gamma]=[\gamma] & \Rightarrow & \mathbbm{1}^0=\textup{Id},\\
\mathbbm{1}^1\big([\gamma]\big)=[(\kappa\cdot e)\ast (e\cdot \gamma)]=[\kappa\ast\gamma] & \Rightarrow & \mathbbm{1}^1=\mathbbm{1}.
\end{eqnarray*}
This implies, in particular, that $\mathbbm{1}_\ast=\textup{Id}$. $\hfill\square$

\subsection{The case of finite fundamental group}\label{ffgp}

\hf\hf Let $n_T(g)$ be the number of geometrically distinct periodic geodesics in $(M,g)$ of length at most $T$. Very few results are known about the asymptotics of the periodic geodesic problem in the case of manifolds with finite fundamental group. In \cite{BTZ81} Ballmann et. al proved the following :
\begin{thm}{\bf (Ballmann, Thorbergsson, Ziller '81)}\hf Let $M$ be a closed, smooth manifold with non-zero but finite fundamental group. Then for any generic smooth metric on $M$ there exists $C>0$ such that
\bgd
n_T\geq C\left(\frac{T}{\log T}\right)
\edd
if $T$ is sufficiently large.
\end{thm}
In particular, the above applies to manifolds of dimension four and for {\it generic} metrics. On the other hand, \cite{VPSu76} tells us that there are {\it infinitely} many geometrically distinct periodic geodesics for {\it any} metric if the manifold has finite fundamental group and the cohomology ring, possibly after passing to a cover, requires at least two generators. We have a result which is in the spirit of these two results in dimension four. \\
\hf\hf We will prove that Corollary \ref{expLM4} and Theorem \ref{expLM} holds for any closed $4$-manifold $M$ with finite fundamental group. As noted before (cf. Corollary \ref{sthomtop2} in \S\ref{sthomgp}), let the second homotopy group be of the form $\pi_2(M)=\Z^k$. 
\begin{cor}\label{expLMffg}
Let $M$ be a closed $4$-manifold with finite fundamental group and $k\geq 3$. For a generic metric the number of geometrically distinct periodic geodesics in $M$ grows exponentially as a function of length.
\end{cor}
{\bf Proof.}\hf Let $\pi:\widetilde{M}\to M$ be the universal cover of $M$. Given a metric $g$ on $M$ we can pull back this metric to a metric $\pi^\ast g$ on the universal cover $\widetilde{M}$. The map $\pi$, being a local isometry, maps closed geodesics to closed geodesics. Moreover, if $\widetilde{\gamma}$ is a closed geodesic in $\widetilde{M}$ and a degenerate critical point of the energy functional then $\pi\circ\widetilde{\gamma}$, a closed geodesic in $M$, is a degenerate critical point of the energy functional. This contradicts the genericity of the metric $g$. Therefore, $\pi^\ast g$ is also a generic metric.\\
\hf Now the covering map $\pi:\widetilde{M}\to M$ is a local isometry and geodesics in $\widetilde{M}$ gets mapped to geodesics in $M$. By the assumption on $k$ we know that geometrically distinct periodic geodesics in $\widetilde{M}$ grow exponentially. However, whenever two geometrically distinct periodic geodesics in the cover maps to the same periodic geodesic in $M$, they are related by Deck transformations. Let $n_T(M,g)$ and $n_T(\widetilde{M},\pi^\ast g)$ be the count of closed geodesics of length at most $T$. Then
\bgd
n_T(M,g)\geq \frac{n_T(\widetilde{M},\pi^\ast g}{|\pi_1(M)|}.
\edd
In particular, $n_T(M,g)$ has exponential growth whenever $k\geq 3$. $\hfill\square$
\begin{rem}
If $M$ is a closed manifold with a fundamental group then $LM$ has connected components indexed by conjugacy classes in $\pi_1(M)$. In particular, we consider $L_0M$, the component containing the constant loops. This is homotopy equivalent to $L\widetilde{M}$, the free loop space of the universal cover. If the Betti numbers of $L\widetilde{M}$ grow exponentially then so does the Betti numbers of $L_0M$ and also of $LM$.
\end{rem}

\subsection{Growth of closed Reeb orbits}\label{Reeb}

\hf\hf The growth of the number of closed geodesics in $M$, arising from the growth of the ranks of homology of the based loop space $\Omega M$, imply analogous results about the growth of Reeb chords for any Reeb flow on a spherization (unit cotangent bundle of $M$).  Similarly, the growth of the number of closed geodesics, arising from the growth of the ranks of homology of the free loop space $LM$, imply results about the growth of closed Reeb orbits. The main technical tool one uses is called {\it sandwiching}. A complete exposition for the case of based loop space (Reeb chords) is given in \cite{MSch11}. For the case of free loop spaces, i.e., closed Reeb orbits, one may refer to either one of two independent sources, work of Macarini, Merry and Paternain \cite{MMPa12}, or the thesis of Heistercamp \cite{Hei11}. In what follows, we provide an extremely brief overview that hopefully explains this connection between closed geodesics and Reeb orbits and then state what we get by combining work done in \cite{Hei11} and \cite{MMPa12} with Theorem \ref{expLM}. We refer the reader to \cite{Hei11}, which is quite readable and is our main source.\\
\hf\hf Let $(M,g)$ be a closed, smooth Riemannian manifold of dimension $m$ and denote by $T^\ast M$ its cotangent bundle. There is a symplectic structure $\omega=d\lambda$, where $\lambda=\sum_{i}^m p_i dq_i$ is the Liouville form written in terms of local coordinates. A hypersurface $\Sigma\subset T^\ast M$ is called {\it fibrewise star-shaped} if for each point $q\in M$ the intersection $\Sigma\cap T^\ast_q M$ is star-shaped with respect to the origin in $T_q^\ast M$. Of particular interest is the example of $\Sigma=T^\ast_1 M$, the unit\footnote{This is one of the instances where we require a Riemannian metric on $M$.} cotangent bundle. Associated to such a hypersurface $\Sigma$, there is vector field $R$ (along $\Sigma$) defined by
\bgd
d\lambda(R,\cdot)=0,\,\,\lambda(R)=1.
\edd
This vector field is called the {\it Reeb vector field} and the flow generated by it is called the {\it Reeb flow}. When $\Sigma=T^\ast_1 M$ the Reeb flow is\footnote{We will identify $TM$ with $T^\ast M$ using the Riemannian metric $g$ and thus a flow on $T^\ast_1 M$ can be compared with the actual geodesic flow on the unit tangent bundle $T_1 M$.} the geodesic flow on $T^\ast_1 M$. For a gentle introduction to geodesic flows and more we refer the reader to \cite{Pat99}. We now note that a metric $g$ being generic (or bumpy) means that the closed geodesics, which are closed orbits of the geodesic flow, are transversally non-degenerate under the geodesic flow. Generalizing this, we call a hypersurface $\Sigma\subset T^\ast M$ {\it generic} if each closed Reeb orbit is transversally non-degenerate under the Reeb flow. When $\Sigma=T^\ast_1 M$ this is equivalent to $g$ being a bumpy metric. Let $\mathpzc{O}_R(\tau)$ be the number of closed Reeb orbits of length at most $\tau$.\\
\hf\hf A fibrewise star-shaped hypersurface $\Sigma\subset T^\ast M$ has a contact structure given by the dual to the hyperplane field
\bgd
\xi_\Sigma:=\ker (\lambda|_\Sigma),
\edd
where $\lambda$ is the Liouville form on $T^\ast M$. It is known that any two such hypersurfaces $(\Sigma,\xi_\Sigma)$ and $(\Sigma',\xi_{\Sigma'})$ are isomorphic as contact manifolds. The identification of these contact manifolds is called the {\it spherization} $(SM,\xi)$ of the cotangent bundle $(T^\ast M,\omega)$. A particular representative, associated to a Riemannian metric $g$ on $M$, is given by $(T^\ast_1 M, \ker \lambda)$. \\
\hf\hf There are several results (cf. \cite{Hei11} Theorem A, B and C) that explain the connections between Reeb orbits and closed geodesics. For our purpose, we shall use the following ( \cite{Hei11} Theorem C, \cite{MMPa12}):
\begin{thm}{\bf (Heistercamp '11; Macarini, Merry, Paternain '12)}\
Suppose that $M$ is a compact, simply connected manifold of dimension at least $2$. Let $\Sigma$ be a generic, fibrewise star-shaped hypersurface of $T^\ast M$ and $R$ its associated Reeb vector field. Then there exists constants $\alpha>0,\beta>0$, depending on $R$, such that
\bgd
\mathpzc{O}_R(\tau)\geq \alpha \max_{i\leq \beta\tau} b_i(LM)
\edd
for $\tau$ sufficiently large.
\end{thm}
This result generalizes \cite{BaZi82} from geodesic flows to Reeb flows; $\Sigma=T^\ast_1 M$ being generic is precisely the case of closed geodesics for a bumpy metric on $M$ and $\mathpzc{O}_R (\tau)=n_\tau(g)$. In view of this theorem and Theorem \ref{expLM} we may conclude the following :
\begin{cor}\label{ReebLM}
Let $M^4$ be a closed, simply-connected, smooth Riemannian manifold with $b_2(M^4)\geq 3$. Let $\Sigma\subset T^\ast M$ a fibrewise star-shaped hypersurface. Then the number of closed Reeb orbits grows at least exponentially as a function of length.
\end{cor}
With $M$ as above, Corollary \ref{ReebLM} may be equivalently stated as "closed Reeb orbits on spherizations of $T^\ast M$ grow at least exponentially for any Reeb flow".


\section{How fast do the homotopy groups grow?}\label{growhg}

\hf\hf In this section we derive an explicit formula for the dimension of rational homotopy groups of simply connected four manifolds in terms of the second Betti number of the manifold. From this formula we can directly deduce that simply connected four manifolds with second Betti number at least $3$ are rationally hyperbolic. The results of this section hold for {\it topological} $4$-manifolds which are not necessarily {\it smooth}. Analogous results were obtained in the case of complete intersections in \cite{Nei79}, \cite{Ba79} using different methods. We shall often use the notation $\pi_k^\Q(X)$ to denote the $k^\text{th}$ rational homotopy group of $X$.


\subsection{Growth of rational homotopy groups}\label{rhg}

\begin{defn}
A path connected topological space $X$ is called {\it nilpotent} if $\pi_1(X)$ is nilpotent and the action of $\pi_1(X)$ on $\pi_n(X)$ is nilpotent.
\end{defn}
Simply connected spaces automatically satisfy the condition above. A space is said to have {\it finite rational cohomology} if 
\bgd
\sum_{i=0}^{\infty}\dim H^i(X;\Q)<\infty.
\edd
\hf\hf One may use rational homotopy theory to divide simply connected manifolds into two classes. 
\begin{defn}
Let $X$ be a nilpotent space with finite dimensional rational cohomology. Then $X$ is called a {\it rationally elliptic space} 
if 
\bgd
\sum_{p\geq 2} \dim \pi_p^\Q(X) <\infty. 
\edd
Otherwise $X$ is called a {\it rationally hyperbolic space}.
\end{defn}
It is well known that rationally elliptic spaces are non-generic and the ellipticity enforces strong conditions (cf. Theorem 2.75 of \cite{RHT208}) on the manifold. In this sense, a generic manifold is rationally hyperbolic. \\
\hf\hf We are interested in understanding the growth rates (of the ranks of the homotopy groups) of rationally hyperbolic spaces among simply connected, closed $4$-manifolds. It is known (cf. Example 3.8, pp. 108-109 of \cite{RHT208}) that the only rationally elliptic, smooth $4$-manifolds are $S^4$, $\CP^2$, $\overline{\CP}^2$, $S^2\times S^2$, $\CP^2\#\CP^2$, and $\CP^2\#\overline{\CP}^2$.\footnote{We have only listed the smooth manifolds in this list. There are rationally elliptic topological $4$-manifolds that are not smoothable but their rational homotopy type agrees with one of the manifolds in this list. In fact, $\CP^2$ and $\overline{\CP}^2$ have the same {\it rational} homotopy type.} Notice that the ellipticity of these manifolds follow very easily from our construction in Lemma \ref{5spin}. In the (smooth) category of simply connected $4$-manifolds, these are precisely the manifolds with second Betti number at most $2$. Consequently, all simply connected $4$-manifolds with second Betti number at least three are rationally hyperbolic. We provide several different proofs of this fact and it's part of a larger circle of ideas.
\begin{prpn}\label{M4kgrow}
Let $M^4_{k}$ be any closed, simply connected $4$-manifold with $H_2(M^4_k;\Z)=\Z^k$. For $i\geq 3$ there is a naturally defined map 
\bgd
\tau:\pi_i(M^4_{k+1})\longrightarrow \pi_i(M^4_{k})
\edd
such that $\tau\otimes\Q$ is a surjection. 
\end{prpn}
{\bf Proof.}\hf To construct the map $\tau$ recall that we have $5$-manifolds $\mathpzc{E}_k$ (respectively $\mathpzc{E}_{k-1}$) which fibre over $M^4_{k+1}$ (respectively $M^4_k$) with fibre $S^1$. Let us denote the bundle maps by $\pi:\mathpzc{E}_k\to M^4_{k+1}$. By construction, $\pi$ induces an isomorphism on homotopy groups in degree $3$ or higher. We also know that $\mathpzc{E}_k$ is diffeomorphic to a connected sum of $k$ copies of $S^2\times S^3$. We have the map
\bgd
q:\mathpzc{E}_k\longrightarrow \mathpzc{E}_{k-1}
\edd
given by the composition of {\it pinching} and then {\it collapsing}. More specifically, let 
\bgd
\textup{pinch}:\mathpzc{E}_k\cong\left(\#^{k-1} S^2\times S^3\right)\# (S^2\times S^3)\longrightarrow \left(\#^{k-1} S^2\times S^3\right)\vee (S^2\times S^3)
\edd
be the map which {\it pinches} the connecting $4$-sphere to a point. Let 
\bgd
\textup{collapse}:\left(\#^{k-1} S^2\times S^3\right)\vee (S^2\times S^3)\longrightarrow \left(\#^{k-1} S^2\times S^3\right)\cong\mathpzc{E}_{k-1}
\edd
be the obvious quotient map which collapses $S^2\times S^3$ to a point. Formally, we define
\bgd
q:=\textup{collapse}\circ\textup{pinch}.
\edd
The map $\tau$ is defined as the composite $\pi_\ast\circ q_\ast\circ (\pi_\ast)^{-1}$. As both $\pi_\ast$'s are isomorphisms (in degree $3$ or higher), $\tau$ is surjective if and only if $q_\ast$ is. \\
\hf The $\textup{collapse}$ map is a retraction onto a wedge summand. Hence it induces a surjection of the homotopy groups. Using minimal models and rational homotopy theory one can prove that (cf. Theorem 3.3 and the proof of Corollary 3.4 in \cite{RHT208}) $\textup{pinch}$ induces a surjection on the rational homotopy groups\footnote{This is not true in general but holds here as our constituent spaces in the connect sum do not have monogenic cohomology rings.}. Combining all of the above, we see that $\tau\otimes\Q$ is a surjection on rational homotopy groups of degree $3$ or higher. $\hfill\square$\\[0.2cm]
\hf\hf We now proceed towards an explicit computation of the group $\pi_j^\Q(M^4_{k})$ for a fixed $k \geq 3$. We start by proving some elementary relations about the generating series of a graded vector space. 
\begin{defn}
Let $V=\oplus_{i\geq 0} V_i$ be a graded vector space (graded by non-negative integers). Define the generating series of $V$ by
\bgd
p_V(t)= \sum_{i\geq 0} \dim (V_i) t^i.
\edd
\end{defn} 
\begin{eg}\label{genExtV}
The free graded commutative algebra $\Lambda V$ on a vector space $V$ is a tensor product of the exterior algebra on the odd dimensional generators and the polynomial algebra on the even dimensional generators. Therefore, whenever $V_0=(0)$ we have
\bgd
p_{\Lambda V}(t)=\frac{\prod_{i\,\textup{odd}} (1+t^i)^{\textup{dim}(V_i)}}{ \prod_{i\,\textup{even}} (1-t^i)^{\textup{dim}(V_i)}}.
\edd
\end{eg}
\begin{eg}\label{genTV}
Let $V$ be a graded vector space with $V_0=(0)$. Let $TV$ denote the tensor algebra of the vector space $V$. The dimension of degree $n$ elements in $TV$ can be calculated as follows. For every ordered partition $(n_1,n_2,\ldots, n_k)$ of $n=\sum_i n_i$, we have a subspace $V_{n_1}\otimes V_{n_2}\otimes \cdots \otimes V_{n_k}$ of $TV$ in grading $n$. Moreover $(TV)_n$ is a direct sum of these terms over distinct ordered partitions of $n$. Let $J(n)$ denote the set of ordered partitions of $n$. Therefore, 
\bgd
\dim (TV)_n = \sum_{(n_1,n_2,\ldots,n_k)\in J(n)} \dim (V_{n_1})\dim (V_{n_2})\cdots \dim (V_{n_k}).
\edd
In terms of generating series we have,
$$p_{TV}(t)= \frac{1}{1-\sum_{i\geq 1} \dim V_i t^i }$$ 
\end{eg}
We have the following proposition whose proof we suspend to \S \ref{details}.
\begin{prpn}\label{geo41}
Let $V$ be the graded vector space with basis $\{x_1,x_2,\ldots,x_k,y_1,y_2,\ldots y_k\}$ where $x_i$'s are in degree $1$ and $y_j$'s are in degree $2$. Let $I\subset T(V)$ be the two sided ideal generated by the element 
\bgd
[x_.,y_.]:=    [x_1,y_1]+[x_2,y_2]+\cdots+[x_k,y_k].
\edd
Then the following holds :\\
\textup{(a)} There is an isomorphism of graded groups
\bgd
H_\lp(\Omega (\#^k S^2\times S^3);\Q)\cong T(V)/I.
\edd
\textup{(b)} The generating series for $T(V)/I$ is given by 
\bgd
 p_{TV/I}(t) = \frac{1}{1-kt-kt^2+t^3}.
\edd
\end{prpn}
\hf\hf Before stating our main computation we recall some facts from rational homotopy theory. Let $X$ be a simply connected topological space. The homotopy groups of the based loop space $\Omega X$ of $X$ are $\pi_i(\Omega X)\cong \pi_{i+1}(X)$. The Whitehead product of homotopy groups of $X$ induces a graded Lie algebra, called the {\it rational homotopy Lie algebra} on $X$, structure on 
\bgd
\pi_\lp^\Q(\Omega X) := \oplus_{i\geq 0} (\pi_i^\Q (\Omega X)).
\edd
For any space $X$, the based loop space $\Omega X$ is a $H$-space. This induces a Pontrjagin product on the homology groups making $H_\lp(\Omega X;\Q):=\oplus_{i\geq 0} H_i(\Omega X;\Q)$ into a ring. Recall the theorem of Milnor and Moore \cite{MM65} :
\begin{thm}\label{MM1}{\bf (Milnor-Moore '65)}\hf For a simply connected space $X$, $H_\lp(\Omega X;\Q)$ is the graded universal enveloping algebra of the graded Lie algebra $\pi_\lp^\Q(\Omega X)$.
\end{thm} 
We are ready to state and prove our result.
\begin{thm}\label{rnkpi4}
Let $M^4_{k}$ be a simply connected, closed manifold of dimension $4$ with $H_2(M^4_{k};\Z)=\Z^{k}$. The rank $m_n(k)$  of $\pi_{n+1}(M^4_{k})\otimes \Q$ is given by
\bge\label{rank}
m_n(k) =\sum_{d|n} (-1)^{n+\frac{n}{d}}\frac{\mu(d)}{d} \sum_{a+2b=n/d}(-1)^b {a+b \choose b}\frac{k^{a}}{a+b} .
\ede
\end{thm}
{\bf Proof.}\hf It suffices to calculate the homotopy groups of $\mathpzc{E}_k=\#^k S^2\times S^3$. The strategy will be to compute the rational homology of the based loop space of $\mathpzc{E}_k$ in two ways and compare them. By Proposition \ref{geo41} we know that
\bgd
H_\lp(\Omega (\#^k S^2\times S^3);\Q)\cong T(V)/I
\edd
where $I,V$ are as in Proposition \ref{geo41}. The ideal $I$ is generated by one element in degree $3$. By Proposition \ref{geo41} the generating series of the graded vector space $T(V)/I$ is 
\bge\label{genfnc1}
p_{TV/I}(t)=\frac{1}{1-kt-kt^2+t^3}. 
\ede
By Theorem \ref{MM1} the rational homology $H_\lp(\Omega \mathpzc{E}_k;\Q)$ is isomorphic (as Hopf algebras) to the universal enveloping algebra of the rational homotopy Lie algebra $L:=\pi^\Q_\lp(\Omega \mathpzc{E}_k)$ of $\Omega \mathpzc{E}_k$. For graded Lie algebras we have a version of Poincar\'{e}-Birkhoff-Witt theorem (cf. Theorem 21.1 of \cite{FHT01}) which says that as an abelian group the universal enveloping algebra $UL$ is isomorphic to $\Lambda L$, the free graded commutative algebra on $L$. Let us denote the dimension of $L_j:=\pi_j^\Q(\Omega \mathpzc{E}_k)$ by $l_j(k)$. Therefore, by Example \ref{genExtV}, $UL$ has the generating series 
\bge\label{genfnc2}
p_{UL}(t)=\frac{\prod_{i\,\textup{odd}} (1+t^i)^{\textup{dim}(L_i)}}{ \prod_{i\,\textup{even}} (1-t^i)^{\textup{dim}(L_i)}}.
\ede
Combining \eqref{genfnc1} with \eqref{genfnc2} we get an equation :
\bgd
\frac{\prod_{i\,\textup{odd}} (1+t^i)^{l_i}}{ \prod_{i\,\textup{even}} (1-t^i)^{l_i}}=\frac{1}{1-kt-kt^2+t^3}=\frac{1}{(1+t)(1-(k+1)t+t^2)}.
\edd
The rank $m_n(k+1)$, also equal to the dimension of $\pi_{n+1}(M^4_{k+1})$, is given by
\bgd
m_n(k+1)=\left\{\begin{array}{rl}
l_n(k) & \textup{if}\,\, n\geq 2\\
l_1(k)+1 & \textup{if}\,\, n=1.
\end{array}\right.
\edd
Therefore, in terms of $m_i(k)$, we have
\bgd
\frac{\prod_{i\,\textup{odd}} (1+t^i)^{m_i(k)}}{ \prod_{i\,\textup{even}} (1-t^i)^{m_i(k)}}=\frac{1}{1-kt+t^2}.
\edd
We take log of both sides :
\begin{eqnarray*}
\log (1-kt+t^2) & = & \sum_{i\,\textup{even}} m_i\log (1-t^i)-\sum_{i\,\textup{odd}} m_i\log (1+t^i)\\
& = & -\sum_{i\,\textup{even}}m_i\left(t^i+\frac{t^{2i}}{2}+\frac{t^{3i}}{3}+\cdots\right)-\sum_{i\,\textup{odd}} m_i\left(t^i-\frac{t^{2i}}{2}+\frac{t^{3i}}{3}+\cdots\right).
\end{eqnarray*}
Expanding this and equating coefficients, we see that 
\bgd
\lambda_n:=\textup{coefficient of $t^n$ in $\log (1-kt+t^2)$} = -\frac{(-1)^n}{n}\bigg(\sum_{d|n} (-1)^d d m_d \bigg).
\edd
We use the M\"{o}bius inversion formula; it gives us 
\bge\label{genfnc4}
m_n =-\sum_{d|n} (-1)^{n+\frac{n}{d}}\mu(d)\frac{\lambda_{n/d}}{d}.  
\ede
In fact, by expanding $\log (1-kt+t^2)$, we see that
\bge\label{lamb}
\lambda_n=-\sum_{a+2b=n}(-1)^b{a+b \choose b}\frac{k^a}{a+b}.
\ede   
This completes the evaluation of $m_n(k)$, the rank of $\pi_{n+1}(M^4_{k})\otimes \Q$.  $\hfill\square$\\[0.2cm]
\hf\hf It follows from the theorem that the ranks of homotopy groups can be recursively computed. We list the first few for the convenience and curiousity of the reader. In particular, it follows very easily that a four manifold as above with second Betti number at least $3$ is rationally hyperbolic. Moreover, this comprehensively improves the main result of Terzi\'{c} \cite{Ter04}, where it is proven that such $4$-manifolds are rationally hyperbolic and explicit ranks of rational homotopy groups (up to degree 4 in the general case) have been computed. 
\begin{cor}\label{4htpygp}
Let $M^4_{k+1}$ be a simply connected, closed $4$-manifold with $H_2(M^4_{k+1};\Z)=\Z^{k+1}$. If $k\geq 2$ then 
\begin{eqnarray*}
\pi_3(M^4_{k+1}) & = & \Z^{\frac{k(k+3)}{2}}\\
\pi_4(M^4_{k+1}) & = & \Z^{\frac{(k-1)(k+1)(k+3)}{3}}\oplus \Z_2^{2k}\\
\pi^\Q_5(M^4_{k+1}) & = & \Q^{\frac{k(k-1)(k+2)(k+3)}{4}}\\
\pi^\Q_6(M^4_{k+1}) & = & \Q^{\frac{k(k-1)(k+1)(k+2)(k+3)}{5}}\\
\pi^\Q_7(M^4_{k+1}) & = & \Q^{\frac{k(k-1)(k+1)(k+3)(k^2+3k+1)}{6}}.
\end{eqnarray*}
\end{cor}
\begin{rem}
The proof of the above is computational. Of interest may be the following :\\
\textup{(1)} The computations of the rank follows from directly computing \eqref{rank}. As for the torsion in the fourth homotopy group, it arises from $\textup{Tor}(\pi_4(\#^k S^2\times S^3))=\Z_2^{2k}$, which can be easily verified via Whitehead products.\\
\textup{(2)} The first five homotopy groups listed may lead one to ask if these polynomials always factor into linear terms. The case of $\pi_7$ answers that negatively.\\
\textup{(3)} It follows from \eqref{genfnc4}, \eqref{lamb} and properties of the M\"{o}bius function that \\
\hf\hf if $n>3$ then $k-1$ divides the rank of $\pi_n^\Q(M^4_{k+1})$;\\
\hf\hf if $n>5$ then $k+1$ divides the rank of $\pi_n^\Q(M^4_{k+1})$;\\
\hf\hf if $n>4$ then $k$ divides the rank of $\pi_n^\Q(M^4_{k+1})$.
\end{rem}
\begin{cor}\label{MM}
For $j\geq 1$ there exists polynomials $P_j(x)$ of degree $j$ and top coefficient $1/j$ such that the rank of $\pi_{j+1}^\Q(M^4_k)$ is $P_{j}(k)$. Moreover,
\bgd
\lim_{j\to \infty} (P_j(k))^{1/j}=\frac{k+\sqrt{k^2-4}}{2}
\edd
whenever $k\geq 3$.
\end{cor}
It follows that $P_n(k)\neq 0$ for $n\geq N$ for some $N$. In particular, $M^{4}_{k+1}$ is rationally hyperbolic whenever $k\geq 2$. \\[0.2cm]
{\bf Proof.}\hf The polynomial $P_j(k)$ can be defined using \eqref{rank}, i.e., let
\bgd
P_n(x):=\sum_{d|n} (-1)^{n+\frac{n}{d}}\frac{\mu(d)}{d} \sum_{a+2b=n/d}(-1)^b {a+b \choose b}\frac{x^{a}}{a+b} .
\edd
Alternatively, $P_n(k)=m_n(k)$ and we may write (refer \eqref{lamb}) $\lambda_n=-\frac{\alpha^n+\beta^n}{n}$, where $\alpha,\beta$ are the roots of the quadratic equation $x^2-kx+1=0$. We also assume that $\beta=(k+\sqrt{k^2-4})/2$ is the bigger root. Note that if $k\geq 3$ then $\alpha/\beta<1$. Therefore,
\bge\label{albe}
nP_n(k)=\sum_{d|n} (-1)^{n+\frac{n}{d}}\mu(d)(\alpha^{n/d}+\beta^{n/d}).
\ede
It then follows that if $k\geq 3$ then $\beta^n$ tends to infinity and
\bgd
\lim_{n\to \infty}\frac{nP_n(k)}{\beta^n}=1.
\edd
This implies, after taking $n^\text{th}$ root of $nP_n(k)/\beta^n$, the limit claimed in the corollary. $\hfill\square$\\[0.2cm]
\hf\hf Recall that there is a dichotomy in action in the fact that rationally hyperbolic spaces have certain growth rates of its (rational) homotopy groups. It is known that if $X$ is a rationally hyperbolic space then the sequence $\sum_{2\leq p\leq m} \dim \pi_p^\Q(X)$ has {\it exponential growth} (cf. Theorem 2.74 of \cite{RHT208}), i.e., there 
are constants $A>1$ and $\lambda>0$ such that, for $m$ large enough, 
\bgd
\sum_{2\leq p\leq m}\dim \pi_p^\Q(X) \geq \lambda A^m.
\edd
\hf\hf With a little work, one can use the recursive formula we have derived to prove the following :
\begin{cor}\label{M432}
Let $M^4_{k+1}$ be a simply connected, closed $4$-manifold with $H_2(M^4_{k+1};\Z)=\Z^{k+1}$. Then for any $n\geq 1$ 
\bgd
\sum_{2\leq j\leq 2n+1}\dim \pi^\Q_j(M^4_{k+1})\geq \frac{k^{2n}}{2n}.
\edd
\end{cor}
{\bf Proof.}\hf We have
\begin{eqnarray*}
\sum_{j=2}^{2n+1} \dim \pi^\Q_j(M^4_{k+1}) & = & \sum_{j=1}^{2n}P_j(k+1)\\
& \geq & \frac{1}{2n}\sum_{j|2n}(-1)^jjP_j(k+1)\\
& = & \frac{1}{2n} (\alpha^{2n}+\beta^{2n})\\
& \geq & \frac{1}{2n}\beta^{2n}\geq \frac{k^{2n}}{2n}.
\end{eqnarray*}
The second equality comes from M\"{o}bius formula applied to \eqref{albe}. Here $\alpha$ and $\beta$ are the roots of $1-(k+1)t+t^2=0$ and $\beta$ is the larger root. The last inequality follows from a simple check as long as $k\geq 1$. $\hfill\square$
\begin{rem}
The fact that $M^4_k$ is rationally hyperbolic for any $k\geq 3$ also follows from Proposition \ref{M4kgrow} by using the fact $M^4_3$ is rationally hyperbolic (cf. Corollary \ref{M432}, or comments following Corollary \ref{MM}). 
\end{rem}


\subsection{Computing the homology of the based loop space}\label{details}

\hf\hf In the last part of this paper we provide a proof of Proposition \ref{geo41}. There are two inputs to this proof, one from rational homotopy theory which reduces the computation of $H_\lp(\Omega (\#^k S^2\times S^3);\Q)$ to an algebraic computation, and the second from Koszul duality of associative algebras which completes the calculation. We begin with a series of definitions. 
\begin{defn}
A (differential) graded Lie algebra $L$ is said to be connected if $L=\oplus_{i\geq 1} L_i$. 
\end{defn}
\begin{defn}
A free model of a connected differential graded Lie algebra $(L,d)$ is a quasi-isomorphism 
\bgd
m:(L(V),d)\rightarrow (L,d)
\edd
such that $L(V)$ is the free graded Lie algebra on the graded vector space $V=\oplus_{i\geq 1}V_i$. 
\end{defn}
\begin{defn}
A Lie model $m:(L(V),d)\rightarrow (L,d)$ is a free model for which $d$ satisfies the minimality condition. The minimality condition states that $d(V)\subset [L(V)^+,L(V)^+]$ where $L(V)^+$ are the positive graded terms of $L(V)$.
\end{defn}
\hf From \cite{FHT01} (cf. Theorem 22.13) we know that $(L(V),d)$ exists for every connected $(L,d)$ and is unique up to isomorphism.
\begin{defn}
Let $X$ be a simply connected space and the rational homology of $X$ is finite dimensional in each degree. Suppose that the rational homotopy Lie algebra is $\mathbb{L}_X$. In this case, a Lie (respectively free) model for $X$ is a Lie (respectively free) model for $(\mathbb{L}_X,0)$.  
\end{defn}
Therefore, for any Lie (free) model $(L(V),d)$ of $X$, the homology $H(L(V),d)$ is the rational homotopy Lie algebra of $X$. When a $(n+1)$-cell is attached to $X$ we have a way to compute the Lie model of the new space in terms of a Lie model of $X$. 
\begin{thm}{\textup{(\cite{FHT01}, Theorem 24.7)}}\label{attach}
Let $Y=X\cup_f \mathcal{D}^{n+1}$ with $f:S^n\rightarrow X$. Suppose $X$ is simply connected and the rational homology is finite dimensional in each degree. Let $(L(V),d)$ be a Lie model for $X$. Let $[f]\in \pi_n(X)=\pi_{n-1}(\Omega X)=H_{n-1}(L(V),d)$ be represented be a cocycle $z$. Then a free model $(L(W),d_W)$ is given by $W=V\oplus \Q\{w_n\}$ and $d_W|_V=d,\, d(w_n)=z$.    
\end{thm}
\hf\hf In our case the above theory proves the following result, which reduces Proposition \ref{geo41} to an algebraic computation. 
\begin{prpn}\label{rathom}
Let $V,I$ be as in Proposition \ref{geo41}. Let $W=V\oplus \Q\{w_4\}$ and define $d:T(W)\rightarrow T(W)$ by $d|_V=0,\, d(w_4)=[x_.,y_.]$. Then
\bgd
H_\lp(\Omega \#^k(S^2\times S^3);\Q) = H_\lp(T(W),d).
\edd
\end{prpn}
{\bf Proof.}\hf The space  $X=\#^k(S^2\times S^3)$ can be given a CW-complex structure with $4$-skeleton the wedge of $k$-copies of $S^2$ ($\alpha_1,\alpha_2,\ldots, \alpha_k$) and $k$-copies of $S^3$ ($\beta_1,\beta_2,\ldots,\beta_k$). There is only one $5$-cell which is attached along the map
\bgd
[\alpha_1,\beta_1]+[\alpha_2,\beta_2]+\ldots+[\alpha_k,\beta_k]
\edd
where the notation $[.,.]$ stands for the Whitehead product. The Lie model can be computed using this CW-complex structure. The Lie model for $\vee^k(S^2\vee S^3)$ is the free Lie algebra  $L(x_1,x_2,\ldots,x_k,y_1,y_2,\ldots,y_k)$ (where $|x_i|=1$, $|y_i|=2$) with 0 differential. Therefore, by Theorem \ref{attach}, the Lie model $(L,d)$ for $X$ is given by $L=L(x_1,x_2,\ldots,x_k,y_1,y_2,\ldots,y_k,w_4)$ with $|w_4|=4$ and the differential generated by
\bgd
d(x_i)=0,~d(y_i)=0,~d(w)=[x_1,y_1]+[x_2,y_2]+\ldots+[x_k,y_k].
\edd
Passing to the universal enveloping algebra we have 
\bgd
H_\lp(\Omega X) = U\mathbb{L}_X = UH(L(W),d) = H(U(L(W)),d) =H(T(W),d)
\edd
as claimed. $\hfill\square$\\[0.2cm]
\hf\hf In order to calculate the above homology we need to use Koszul duality of associative algebras (cf. \cite{LoVa12}). 
\begin{defn}
The associative algebra $A(V,R)=T(V)/(R)$ is a {\it quadratic} algebra if $R$ can be chosen as $R\subseteq V\otimes V$. 
\end{defn} 
The dual notion leads to a quadratic coalgebra. For this note that $T(V)$ has a coalgebra structure by declaring the elements of $V$ to be primitive\footnote{An element $c$ in a coalgebra $ C$ is called primitive if $\Delta(c)=c\otimes 1 + 1\otimes c$.} . When $T(V)$ is thought of as a coalgebra we write it as $T^c(V)$. 
\begin{defn}
For $R\subset V\otimes V$ the quadratic coalgebra $C(V,R)$ is defined as the universal object in the category of sub-coalgebras $C$ of $T^c(V)$ so that $C\to T^c(V) \to V\otimes V/R$ is $0$. The universal property of $C(V,R)$ states that for any sub-coalgebra $C$ of $T^c(V)$ so that the weight $2$ elements are contained in $R$, $C\subset C(V,R)$. 
\end{defn}
\hf For an augmented graded algebra $A$, there is a bar construction $B(A)$ which is a differential graded coalgebra and dually for a coaugmented graded coalgebra $C$ there is a cobar construction $\Omega C$ which is a differential graded algebra. These are defined as   
\begin{defn}
Let $\bar{A}\subset A$ be the kernel of the augmentation, define $B(A)= (T(s\bar{A}),d)$ ($s=$ suspension) where $d$ is generated as a coderivation by
\bgd
d(s(a))=s(a\otimes a)-s(a\otimes 1) -s(1\otimes a).
\edd
Dually, let $C=\bar{C}\oplus \Q$, and define $\Omega C = (T(s^{-1}\bar{C}),d)$ where $d$ is generated as a derivation by the equation 
\bgd
d(s^{-1} c) = s^{-1}(\bar{\Delta} (c))=s^{-1}(\Delta(c)-c\otimes 1 - 1\otimes c).
\edd
\end{defn}
\hf\hf A quadratic algebra (respectively coalgebra) as defined is augmented and the above constructions apply. There is a differential on $C\otimes \Omega C$ generated by $d(c)=1\otimes s^{-1}c$ and dually a differential on $A\otimes B(A)$.  
\begin{defn}
The Koszul dual coalgebra of a quadratic algebra $A(V,R)$ is defined as $A^\textsuperscript{!`}=C(s(V),s^2(R))$. The Koszul dual algebra $A^\textup{!}$ of a quadratic algebra $A(V,R)$ is defined as $A^\textup{!}=A(V^*,R^\perp)$ where $R^\perp\subset V^*\otimes V^*$ consists of elements which take value $0$ on $R\subset V\otimes V$. 
\end{defn}
The Koszul dual algebra and the Koszul dual coalgebra are linear dual up to a suspension. Let $A^{(n)}$ stand for the subspace of homogeneous $n$-fold products. Then $(A^!)^{(n)}\cong s^n((A^\textsuperscript{!`})^*)^{(n)}$. \\
\hf\hf For a quadratic algebra $A(V,R)$, there is a natural map from $\Omega A^\textsuperscript{!`} \rightarrow A$ which maps $v\mapsto v$. Using this map there is a differential on $A^\textsuperscript{!`}\otimes A$ denoted by $d_\kappa$. We shall use the notation $A^\textsuperscript{!`}\otimes_\kappa A$ to denote $(A^\textsuperscript{!`}\otimes A,d_\kappa)$. 
\begin{defn}{(cf. Theorem 3.4.6 of \cite{LoVa12})}\label{Kos}
A quadratic algebra $A(V,R)$ is called Koszul if $\Omega A^\textsuperscript{!`}\rightarrow A$ is a quasi-isomorphism. Equivalently, the chain complex $A^\textsuperscript{!`}\otimes_\kappa A$ is acyclic. 
\end{defn}
\hf\hf There is an easy way to check that an algebra is Koszul. First fix a basis $(v_1,v_2,\ldots, v_m)$ of $V$, and fix an order $v_1<v_2\cdots<v_m$. This induces a lexicographic order on the degree 2 monomials. Now arrange the expressions in $R=(r_1,r_2,\ldots)$ in terms of order of monomials. An element $v_iv_j$ is called a leading monomial if there exists $r_l=v_iv_j + \mbox{lower order terms}$. Note that (\cite{LoVa12}, Theorem 4.1.1) implies that if there is only one leading monomial $v_iv_j$ with $i\neq j$ then the algebra is Koszul. We use this to prove the following :
\begin{prpn}\label{Kos2}
Let $V,I$ be as in Proposition \ref{geo41}. Then the algebra $A=A(V,I)=T(V)/I$ is Koszul.
\end{prpn} 
{\bf Proof.}\hf We can order the monomials as $x_1<x_2\ldots <x_k< y_1<\ldots < y_k$. Then there is only one leading monomial $y_nx_n$. This completes the proof by the above criterion. $\hfill\square$\\[0.2cm]
\hf\hf Let $A=A(V,I)$ as in Proposition \ref{geo41}. The Koszul dual algebra of $A$ is the algebra $T(V^*)/(I^\perp)$, where $I^\perp\subset V^*\otimes V^*$ is the subset which acts by $0$ on the element $[x_.,y_.]$. This is spanned by 
\bgd
\{x_i\otimes x_j,~y_i\otimes y_j,~x_i\otimes y_j (j\neq i),~x_i\otimes y_i-x_j\otimes y_j,~x_i\otimes y_i-y_i\otimes x_i\}.
\edd
Therefore, the Koszul dual algebra of $T(V)/I$ is $\Q \oplus V^* \oplus \Q\{r^*\}$ where 
\bgd
x_i^*x_j^*=0,~y_i^*y_j^*=0,~x_i^*y_j^*=y_j^*x_i^*=\delta_{ij} r^*.
\edd
\begin{prpn}\label{cobar}
Let $V,I$ as in Proposition \ref{geo41}. Let $(C,\Delta)$ be the coalgebra : $\Q \oplus s(V)\oplus \Q\{s(w)\}$ with 
\begin{eqnarray*}
\Delta(s(v)) & = & s(v)\otimes 1 + 1\otimes s(v)~\mbox{for}~ v\in V\\
\Delta(s(w)) & = & s(w)\otimes 1 + 1\otimes s(w) + [s(x)_.,s(y)_.].
\end{eqnarray*}
Then the cobar complex $\Omega C$ equals the differential graded algebra $(T(W),d)$ of Proposition \ref{rathom}.
\end{prpn}
{\bf Proof.}\hf The cobar construction $\Omega C$ as an algebra is given by $T(V\oplus \Q\{w\})$ which is isomorphic to $W$. Since $s(v)$'s are primitive, $d(v)=0$ and
\bgd
d(w) = s^{-1}\bar{\Delta}(w)= s^{-1}([sx_.,sy_.])= [x_.,y_.].
\edd
$\vspace*{-0.8cm}\hfill\square$\\[0.9cm]
Finally, the following proposition follows easily from the computation of the Koszul dual algebra of $A$.
\begin{prpn}\label{coalg}
The coalgebra $C$ of Proposition \ref{cobar} is the Koszul dual coalgebra of the algebra $A=A(V,I)$. 
\end{prpn}
We're now ready to prove Proposition \ref{geo41}.\\[0.2cm]
{\bf Proof of Proposition \ref{geo41} (a) :} By Proposition \ref{rathom}, we have
\bgd
H_\lp(\Omega \#^k(S^2\times S^3);\Q) = H_\lp(T(W),d)
\edd
where $W=V\oplus \Q\{w_4\}$ and $d|_V=0,\, d(w_4)=[x_.,y_.]$. We need to prove that this homology is $A:=T(V)/I$. By Proposition \ref{Kos2}, the algebra $A$ is Koszul. By proposition \ref{coalg} the Koszul dual coalgebra of $A$ is the coalgebra $C$ of Proposition \ref{cobar}. This implies that the map $\Omega C \rightarrow A$ is a quasi-isomorphism. By Proposition \ref{cobar}, the cobar construction $\Omega C$ is isomorphic to the differential graded algebra $(T(W),d)$. Therefore, we have that $H_\lp(T(W),d) \cong A = T(V)/I$ as required. $\hfill \square$ \\[0.2cm]
{\bf Proof of Proposition \ref{geo41} (b) :} Let $C$ be as in Proposition \ref{cobar}. Proposition \ref{Kos2} implies that $A$ is Koszul and Proposition \ref{cobar} implies that $C$ is its Koszul dual coalgebra. Hence, the Koszul complex $C\otimes_\kappa A$ is acyclic. The linear map $\kappa$ is defined as
\bgd
\kappa:C\longrightarrow A,\,\,s(v)\mapsto v,\,s(w)\mapsto 0,\,1\mapsto 0.
\edd
The differential $d_\kappa$ is the $A$-derivation generated by 
\bgd
C\xrightarrow{\Delta} C\otimes C \xrightarrow{\textup{Id} \otimes \kappa}C\otimes A.
\edd
We expand this out completely using $C=\Q\oplus s(V)\oplus \Q\{s(w)\}$. The vector space $V$ is graded with exactly two non-zero degrees $1$,$2$ and $\dim(V_1)=\dim(V_2)=k$. Let $A_i$ be the $i^\text{th}$ graded part of $A$. Now observe that
\begin{eqnarray*}
d_\kappa : \Q\otimes A_k & \rightarrow & 0\\
d_\kappa : s(V_1)\otimes A_k & \rightarrow & A_{k+1}\\
d_\kappa : s(V_2)\otimes A_k & \rightarrow & A_{k+2}\\
d_\kappa : \Q\{s(w)\}\otimes A_k & \rightarrow & s(V_1)\otimes A_{k+2}\oplus s(V_2)\otimes A_{k+1}.
\end{eqnarray*}
\hf We consider the subcomplexes 
\bgd
C_k=\Q\otimes A_k \oplus s(V_1)\otimes A_{k-1}\oplus s(V_2) \otimes A_{k-2}\oplus \Q\{s(w)\}\otimes A_{k-3}
\edd 
of the chain complex $C\otimes_\kappa A$ which is acyclic. This means that the homology of $C_k$ is $0$ if $k>0$ and $\Q$ if $k=0$. Therefore the Euler characteristic of $C_k$ is
\bgd
\dim(A_k)-k\dim(A_{k-1})-k\dim(A_{k-2})+\dim(A_{k-3})=\left \{\begin{array}{rl} 
                                      0 &\mbox{if}~k>0 \\
                                      1 &\mbox{if}~k=0.  
                                       \end{array}\right.
\edd
\hf In terms of generating series this means
\bgd 
p_A(t)(1-kt-kt^2+t^3)=1
\edd
and consequently,
\bgd
p_A(t) = \frac{1}{1-kt-kt^2+t^3}.
\edd
$\hfill \square$


\bibliographystyle{siam}

\vspace*{0.4cm}

\hf {\small D}{\scriptsize EPARTMENT OF }{\small M}{\scriptsize ATHEMATICS, }{\small R}{\scriptsize AMAKRISHNA }{\small M}{\scriptsize ISSION }{\small V}{\scriptsize IVEKANANDA }{\small U}{\scriptsize NIVERSITY, }{\small H}{\scriptsize OWRAH, }{\small WB} {\footnotesize 711202, }{\small INDIA}\\
\hf{\it E-mail address} : \texttt{samik.basu2@gmail.com}\\[0.2cm]

\hf {\small D}{\scriptsize EPARTMENT OF }{\small M}{\scriptsize ATHEMATICAL }{\small S}{\scriptsize CIENCES, }{\small B}{\scriptsize INGHAMTON }{\small U}{\scriptsize NIVERSITY, }{\small B}{\scriptsize INGHAMTON, }{\small NY} {\footnotesize 13902-6000, }{\small USA}\\
\hf{\it E-mail address} : \texttt{somnath@math.binghamton.edu}

\end{document}